\documentclass[12pt]{article}
\usepackage{latexsym}
\usepackage{amsmath}
\usepackage{amssymb}
\usepackage{amscd}
\usepackage{graphicx}
\begin{document}

\newtheorem{Th}{Theorem}[section]
\newtheorem{Cor}{Corollary}[section]
\newtheorem{Prop}{Proposition}[section]
\newtheorem{Lem}{Lemma}[section]
\newtheorem{Def}{Definition}[section]
\newtheorem{Rem}{Remark}[section]
\newtheorem{Ex}{Example}[section]
\newtheorem{stw}{Proposition}[section]


\newcommand{\bet}{\begin{Th}}
\newcommand{\ent}{\stepcounter{Cor}
   \stepcounter{Prop}\stepcounter{Lem}\stepcounter{Def}
   \stepcounter{Rem}\stepcounter{Ex}\end{Th}}


\newcommand{\bec}{\begin{Cor}}
\newcommand{\enc}{\stepcounter{Th}
   \stepcounter{Prop}\stepcounter{Lem}\stepcounter{Def}
   \stepcounter{Rem}\stepcounter{Ex}\end{Cor}}
\newcommand{\bep}{\begin{Prop}}
\newcommand{\enp}{\stepcounter{Th}
   \stepcounter{Cor}\stepcounter{Lem}\stepcounter{Def}
   \stepcounter{Rem}\stepcounter{Ex}\end{Prop}}
\newcommand{\bel}{\begin{Lem}}
\newcommand{\enl}{\stepcounter{Th}
   \stepcounter{Cor}\stepcounter{Prop}\stepcounter{Def}
   \stepcounter{Rem}\stepcounter{Ex}\end{Lem}}
\newcommand{\bef}{\begin{Def}}
\newcommand{\enf}{\stepcounter{Th}
   \stepcounter{Cor}\stepcounter{Prop}\stepcounter{Lem}
   \stepcounter{Rem}\stepcounter{Ex}\end{Def}}
\newcommand{\ber}{\begin{Rem}}
\newcommand{\enr}{
   \stepcounter{Th}\stepcounter{Cor}\stepcounter{Prop}
   \stepcounter{Lem}\stepcounter{Def}\stepcounter{Ex}\end{Rem}}
\newcommand{\bee}{\begin{Ex}}
\newcommand{\ene}{
   \stepcounter{Th}\stepcounter{Cor}\stepcounter{Prop}
   \stepcounter{Lem}\stepcounter{Def}\stepcounter{Rem}\end{Ex}}
\newcommand{\Proof}{\noindent{\it Proof\,}:\ }

\renewcommand{\normalbaselines}{\baselineskip20pt
   \lieskip3pt \lineskiplimit3pt}
\newcommand{\mapupright}[1]{\smash{\mathop{
   \hbox to 1cm{\rightarrowfill}}\limits^{#1}}}
\newcommand{\mapdownright}[1]{\smash{\mathop{
   \hbox to 1cm{\rightarrowfill}}\limits_{#1}}}
\newcommand{\mapdown}[1]{\Big\downarrow 
   \llap{$\vcenter{\hbox{$\scriptstyle#1\,$}}$ }} 
\newcommand{\mapup}[1]{\Big\uparrow 
   \rlap{$\vcenter{\hbox{$\scriptstyle#1\,$}}$ }} 
\newcommand{\EE}{\mathcal{E}}
\newcommand{\KK}{\mathbf{K}}
\newcommand{\QQ}{\mathbf{Q}}
\newcommand{\R}{\mathbf{R}}
\newcommand{\C}{\mathbf{C}}
\newcommand{\ZZ}{\mathbf{Z}}
\newcommand{\NN}{\mathbf{N}}
\newcommand{\PP}{\mathbf{P}}
\newcommand{\OOO}{{\mathbf{O}}}
\newcommand{\HHH}{{\mathbf{H}}}
\newcommand{\AAA}{{\mathbf{A}}}
\newcommand{\uuu}{\boldsymbol{u}}
\newcommand{\xxx}{\boldsymbol{x}}
\newcommand{\aaa}{\boldsymbol{a}}
\newcommand{\bbb}{\boldsymbol{b}}
\newcommand{\BBB}{\mathbf{B}}
\newcommand{\ccc}{\boldsymbol{c}}
\newcommand{\iii}{\boldsymbol{i}}
\newcommand{\jjj}{\boldsymbol{j}}
\newcommand{\kkk}{\boldsymbol{k}}
\newcommand{\rrr}{\boldsymbol{r}}
\newcommand{\FFF}{\mathbf{F}}
\newcommand{\yyy}{\boldsymbol{y}}
\newcommand{\ppp}{\boldsymbol{p}}
\newcommand{\qqq}{\boldsymbol{q}}
\newcommand{\nnn}{\boldsymbol{n}}
\newcommand{\vvv}{\boldsymbol{v}}
\newcommand{\eee}{\boldsymbol{e}}
\newcommand{\fff}{\boldsymbol{f}}
\newcommand{\www}{\boldsymbol{w}}
\newcommand{\0}{\boldsymbol{0}}
\newcommand{\lon}{\longrightarrow}
\newcommand{\ga}{\gamma}
\newcommand{\pa}{\partial}
\newcommand{\QED}{\hfill $\Box$}
\newcommand{\id}{{\mbox {\rm id}}}
\newcommand{\Ker}{{\mbox {\rm Ker}}}
\newcommand{\grad}{{\mbox {\rm grad}}}
\newcommand{\ind}{{\mbox {\rm ind}}}
\newcommand{\Gr}{{\mbox {\rm Gr}}}
\newcommand{\rank}{{\mbox {\rm rank}}}
\newcommand{\sign}{{\mbox {\rm sign}}}
\newcommand{\Symp}{{\mbox {\rm Sp}}}
\newcommand{\weight}{{\mbox {\rm weight}}}
\newcommand{\ord}{{\mbox {\rm ord}}}
\newcommand{\GG}{{\mathcal{G}}}
\newcommand{\un}{{\mathrm{un}}}
\newcommand{\FF}{{\mathcal{F}}}
\newcommand{\Int}{{\mbox {\rm Int}}}
\newcommand{\codim}{{\mbox {\rm codim}}}
\newcommand{\rest}{{\mbox{\rm \tiny{rest.}}}}
\def\mod{{\mbox {\rm mod}}}
\newcommand{\enP}{\hfill $\Box$ \par\vspace{5truemm}}
\newcommand{\qed}{\hfill $\Box$ \par}

\title{
Infinitesimal deformations and stabilities 
of \\ 
singular Legendre submanifolds. 
}

\author{Go-o ISHIKAWA}
\date{ }
\renewcommand{\thefootnote}{\fnsymbol{footnote}}
\maketitle

\section{Introduction.}  
\label{Intro}

\begin{abstract}
We give the characterization 
of Arnol'd-Mather type for 
stable singular Legendre immersions. 
The most important building block of the theory 
is providing a module structure 
on the space of infinitesimal integral deformations by 
means of the notion of natural liftings of differential systems 
and of contact Hamiltonian vector fields. 
\end{abstract}

The framework of Legendre singularity theory 
for Legendre immersions is established 
\cite{AGV}: 
The singularity of  a Legendre immersion via a Legendre fibration 
is embodied in a family of hypersurfaces, namely, 
the generating family of the Legendre immersion, 
and the stability of such singularity is expressed by 
mean of a notion, ${\mathcal K}$-versality, for its generating family. 
However, since a singular Legendre immersion has no generating family in general, 
the direct characterization 
should be worthwhile for the understanding of the 
stability of singular Legendre immersions, 
which we are going to provide in this paper. 

\

The significance of Legendre singularity theory has increased 
recently by the trend of differential geometry treating (wave) fronts 
as generalised objects of hypersurfaces. Moreover the point of view  
in the micro-local analysis provides the motivation for the study of 
Legendre submanifolds as the description of singularities of solutions 
to partial differential equations. 

The most simple singularity of front is given by 
$$
x = t^2, \quad y = t^3, 
$$
near $(x, y) = (0, 0)$, 
the $(2, 3)$-cusp on the $(x, y)$-plane. 
The front lifts to the Legendre curve
$$
x = t^2, \quad y = t^3, \quad p = \dfrac{3}{2}t, 
$$
which is an immersion. 
Then the stability of the front is well described by the lifted non-singular 
Legendre submanifold via the Legendre equivalences induced 
by diffeomorphisms on the $(x, y)$-plane. 

Consider then the similarly simple plane curve 
$$
x = t^2, \quad y = t^5, 
$$
the $(2, 5)$-cusp near $(x, y) = (0, 0)$. 
Then the natural lifting has the form: 
$$
x = t^2, \quad y = t^5, \quad p = \dfrac{5}{2}t^3, 
$$
which is an integral curve to the contact distribution 
$dy - pdx = 0$ and not an immersion at $t = 0$. 
Therefore, restricted ourselves to Legendre immersion without singularities 
we can not treat this very simple curve in the framework of 
Legendre singularity theory. 
Thus we are going to study, 
in this paper, singular Legendre immersions, in particular, 
the nature of their deformations in canonical way. 

\begin{figure}[htbp]
  \begin{center}
      \includegraphics[width=8truecm, clip, keepaspectratio]{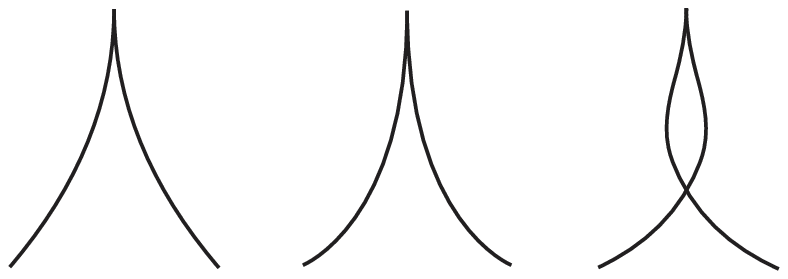} 
    \caption{An integral deformation of $(2, 5)$-cusp.}
    \label{integral}
  \end{center}
\end{figure}

For example, consider the \lq\lq stable\rq\rq deformation of the $(2, 5)$-cusp 
$$
x = t^2, \quad y = t^5 + \lambda t^3
$$ 
inducing smooth deformation 
of tangent lines. See figure \ref{integral}. 

Then we understand, via our general theory, the stable deformation forms the 
stable projection (front) of the open Whitney umbrella of type $1$, introduced 
in this paper, which is contactomorphic to   
$$
x = t^2, \quad y = t^5 + \lambda t^3, \quad p = \dfrac{5}{2}t^3 + \dfrac{3}{2}\lambda t, 
\quad \mu = t^3, 
$$
in the $(x, y, \lambda, p, \mu)$-space with the contact 
structure $dy - pdx - \mu d\lambda = 0$. 

\

We have given in \cite{Inv} 
the characterisations for symplectic stability of Lagrange varieties 
and Lagrange stability in 
symplectic geometry. 
Therefore the present paper can be regarded as a contact or Legendre 
counterpart 
to \cite{Inv}: 
We observe surely the parallelism between Lagrange 
and Legendre singularity theories, as well as symplectic and contact 
geometries. In fact we use several results in Lagrange singularities proved 
in \cite{Inv} to deduce several results in Legendre singularities. 
Nevertheless we need to break through several difficulties 
for obtaining the characterisations (Theorem \ref{Stability}). 

In particular, we realize that 
the direct characterisation needs the deep understanding 
of {\it the space of Legendre submanifolds}. 
Since the space of submanifolds can be treated as the space of immersions, 
we consider, in a contact manifold, the space of integral mappings, 
parametrizations of 
integral submainfolds of the contact distribution. 
The space of integral mappings turns out to be our central object. 
Its tangent space at an integral mapping is naturally regarded 
as the space of infinitesimal deformations of the integral mapping 
among integral mappings. 
The fact, then, we observe in this paper is that the tangent space 
to the space of integral mappings has the structure of not merely a 
vector space, but the very natural 
module structure. 
It reminds us the
 \lq\lq modularity\rq\rq in the sense of Mather \cite{Mather2}. 
However, in this paper we introduce the module structure for 
functions not on the source manifold but for functions on the target manifold. 

We understand the modularity of tangent spaces to the space of integral 
mappings in a contact manifold without difficulty as follows: 
An infinitesimal deformations on a 
contact manifold, namely,  
a contact vector field is locally given by a 
{\it contact Hamilton vector field} $X_K$ with a Hamiltonian  
function $K$ on the contact manifold, fixing a local contact form 
$\alpha$. Then we see, for functions $H, K$,  
$$
X_{HK} = H\cdot X_K + K\cdot X_H - (HK)\cdot X_1. 
$$
Thus, we can give the module structure $*$ of functions on 
the space of contact Hamilton vector fields by identifying it 
with the space of functions: The formula reads as 
$$
H*X_{K} := H\cdot X_K + K(X_H - H\cdot X_1). 
$$
Note that the interior product $i_{X_K}\alpha$ is equal to $K$. 
Let $f$ be an integral mapping. 
The vector field $X_K\circ f$ along $f$ 
is a kind of integral infinitesimal deformations of $f$. 
Then we set 
$$
H*(X_{K}\circ f) := f^*H\cdot (X_K\circ f) + f^*K(X_H - H\cdot X_1)\circ f. 
$$
Note that $i_{X_K\circ f}\alpha$ is equal to $f^*K$. 
We proceed even further. 
We define 
the multiplication $*$ by a function $H$ 
of any infinitesimal integral deformation $v$ 
of $f$ by the formula: 
$$
H*v := f^*H\cdot v + (i_v\alpha)(X_H - H\cdot H_1)\circ f. 
$$
Moreover we observe that the multiplication is {\it intrinsic}: The 
definition of multiplications 
looks like depending on the choice of a local contact form $\alpha$, 
but in fact it is independent of it and is determined only by the 
contact structure. 

Note that the module structure is effectively used in \cite{CSLCC} for 
the classifying of singular Legendre curves in the contact three space. 

We introduce the class of singular Legendre submanifolds, 
{\it open Whitney umbrellas},  
in contact manifolds by explicit forms, and formulate the 
characterisations of Legendre stability and Legendre 
versality; the main Theorems \ref{Stability} and \ref{Versality} in \S \ref{Main results.}. 
In \S \ref{contact stability}, we give the characterisation of 
open Whitney umbrellas as contact stable integral map-germs of corank at most one. 
To prove Theorem \ref{Stability}, 
we need to 
clarify the infinitesimal condition on Legendre stability. 
For this, we introduce the notion of 
natural liftings (\cite{YI}\cite{YK}) of differential forms 
and differential systems 
 in \S \ref{Lie derivative.}. 
After reviewing the notion of contact Hamilton 
vector fields in \S \ref{Contact Hamilton vector fields.}, 
we formulate exactly infinitesimal conditions in \S \ref{Infinitesimal deformations.}. 
In \S \ref{Isotropic map}, 
we study the relation of integral mappings and isotropic mappings, 
and,  
in \S \ref{Integral jet}, we study on the integral jet spaces. 
In \S \ref{Finite determinacy}, we give results on finite determinacy 
of integral map-germs. 
We give, using all results given in the previous sections, 
the proof of Theorem \ref{Stability} in \S \ref{Legendre}. 
In \S \ref{versality}, we mention on the proof of Legendre versality theorem 
\ref{Versality}. 

In this paper, 
all manifolds and mappings we treat are assumed to be of class 
$C^\infty$ (in case $\KK = \R$) or complex analytic 
(in case $\KK = \C$).

\ 

This paper is based on my talk at the workshop in the University of Warwick 
held on June 1999. I would like to thank the organisers Professor 
David Mond and Professor Andrew du Plessis. 
Also I would like to thank Professor Jim Damon 
for valuable comment on the method described in the last section given 
to me at Newton Institute, Cambridge, 2000. 
I would like to thank Professor Hajime Sato for helpful comment, in particular, 
for his indicating me 
the notion of the natural liftings in \cite{YI}. 
Also I would like to dedicate the present paper to Professor Syuzo Izumi 
for his 65th birthday and Professor Hajime Sato for his 60th birthday. 

\section{Main results.}
\label{Main results.}

Now we are going to describe in detail the objects we apply our theory, 
before formulating the main Theorem \ref{Stability}. 

Let $(W, D)$ be a real or complex contact manifold of dimension $2n + 1$ 
\cite{A72}\cite{A90}\cite{AGV}. 
Here $D \subset TW$ stands for the contact 
structure on $W$, namely, a completely 
non-integrable distribution of codimension one.  
A typical example is $W = \KK^{2n+1}$, $\KK = 
\R$ or $\C$, with 
coordinates $(p, q, r)$, and 
$$
D = \{ dr - \sum_{i=1}^n p_idq_i = 0\} \subset T\KK^{2n+1}. 
$$
By Darboux's theorem, any contact manifold is locally contactomorphic to 
this standard model. 

A mapping $f : N \to W$ from an $n$-dimensional manifold $N$ is called 
an {\it integral mapping} if, for any $x \in N$, 
$f_*(T_xN) \subset D_{f(x)}$, where $f_* : T_xN \to T_{f(x)}W$ is 
the differential mapping (the linearization) 
of $f$ at $x$. Thus 
the notion of integral mappings generalizes that of 
(immersed) integral manifolds in the contact manifold $W$.

Two map-germs $f : (N, x_0) \to (W, D)$ and $f' : (N', x'_0) \to (W', D')$ 
to contact manifolds $(W, D)$ and $(W', D')$ respectively, 
are called 
{\it contactomorphic} 
if there exist a diffeomorphism $\sigma : (N, x_0) \to (N', x'_0)$ 
and a contactomorphism $\tau : (W, f(x_0)) \to 
(W', f'(x'_0))$, $\tau_*D = D'$, 
such that $f'\circ\sigma = \tau\circ f$. 
In this case we call also the pair $(\sigma, \tau)$ 
a {\it contactomorphism} of $f$ and $f'$. 

Let $f : (N^n, x_0) \to W^{2n+1}$ be an integral map-germ. 
Suppose that $f$ is of corank $\leq 1$, namely that 
the kernel of the differential map $f_* : T_{x_0}N 
\to T_{f(x_0)}W$ is zero or one dimensional. 
Then there exists a contactomorphism $(\sigma, \tau)$ from $f$ 
to $f' = \tau\circ f \circ\sigma : (\KK^n, 0) \to (\KK^{2n+1}, 0)$ 
such that 
$$
(q_1, \dots, q_{n-1}, q_n)\circ f' = (x_1, \dots, x_{n-1}, u(x_1, \dots, x_{n-1}, x_n), 
$$
for some function $u$, 
where $(x_1, x_{n-1}, x_n)$ is the standard coordinate of $\KK^n$. 
Then, setting $v := p_n\circ f'$, we easily see 
that the components $p_1, \dots, p_{n-1}$ and $r$ of $f'$ are 
uniquely determined by the condition
$$
d(r\circ f) = \sum_{i=1}^n(p_i\circ f)d(q_i\circ f). 
$$
Actually we have: 
\bep 
\label{pre-normal} 
{\rm (Pre-normal form of integral map-germ of corank at most one.)} 
Let $f : (N^n, x_0) \to W^{2n+1}$ be an integral map-germ of corank $\leq 1$. 
Then there exist functions-germs $u, v : (\KK^{n}, 0) \to (\KK, 0)$ 
such that $f$ is contactomorphic to the integral map-germ 
$f' : (\KK^n, 0) \to (\KK^{2n+1}, 0)$ defined by 
$$
(q_1, \dots, q_{n-1}, q_n, p_n)\circ f' := (x_1, \dots, x_{n-1}, u, v), 
$$
$$
p_i\circ f' := \int_0^{x_n} \left( \dfrac{\pa v}{\pa x_i}\dfrac{\pa u}{\pa x_n} - 
\dfrac{\pa v}{\pa x_n}\dfrac{\pa u}{\pa x_i}\right) dx_n, (1 \leq i \leq n-1), 
$$
and
$$
r\circ f' := \int_0^{x_n} \left( v\dfrac{\pa u}{\pa x_n}\right) dx_n. 
$$
\enp

In particular, our main objects of the study are introduced as follows: 
For an integer $k$ with $0 \leq k \leq \dfrac{n}{2}$, we 
define a map-germ $f = f_{n,k} : (\KK^n,0) \to 
(\KK^{2n+1}, 0)$ by 
$q_1\circ f = x_1, \dots, q_{n-1}\circ f = x_{n-1}$ and 
\[ 
u = q_n\circ f = \frac{x_n^{k+1}}{(k+1)!} + x_1\frac{x_n^{k-1}}{(k-1)!} + 
\cdots  + x_{k-1}x_n, 
\]
\[
v = p_n\circ f = x_k\frac{x_n^k}{k!} + \cdots  + x_{2k-1}x_n, 
\]
and the property $f^*\alpha = 0$. 
The components $p_1, \dots, p_{n-1}$ and $r$ of $f$ are 
defined as in Proposition \ref{pre-normal} so that 
$$
d(r\circ f) = \sum_{i=1}^n(p_i\circ f)d(q_i\circ f). 
$$
Then we call a map-germ 
$f : (N, x_0) \to W$ an 
{\it open Whitney umbrella} (or an {\it unfurled Whitney umbrella}) 
{\it of type} $k$ ($0 
\leq k \leq \dfrac{n}{2}$), 
if it is contactomorphic to the normal form $f_{n, k}$. 

An open Whitney umbrella is an integral map-germ of corank at most one.  
It is an immersion, namely, {\it Legendre immersion},  exactly when $k = 0$: 
A map-germ $f : (N^n, x_0) \to (W^{2n+1}, D)$ 
is a Legendre immersion if and only if $f$ is an open Whitney umbrella 
of type $0$. 
If $k > 0$, then the singular locus of 
an open Whitney umbrella of type $k$ is 
non-singular and of codimension $2$ in $N$. 

Open Whitney umbrellas are intrinsically characterised 
via the notion of \lq\lq contact stability" in \S \ref{contact stability}. 

\smallskip  

A fibration $\pi : W^{2n+1} \to Z^{n+1}$ 
is called a {\it Legendre fibration} 
if the fibers of $\pi$ are Legendre submanifolds of $W$. 
Then 
we concern with the relative position of 
the image of an integral mapping  
with respect to a Legendre 
fibration: 
We consider an integral map-germ $f : (N, x_0) 
\to (W, w_0)$ together with a germ of Legendre 
fibration $\pi : (W, w_0) \to (Z, z_0)$, 
where $w_0 = f(x_0), z_0 = \pi(w_0) = (\pi\circ f)(x_0)$. 

Let $\pi : (W, w_0) \to (Z, z_0)$ and 
$\pi' : (W', w'_0) \to (Z', z'_0)$ be germs of Legendre fibrations. 
Then 
a contactomorphism-germ $\tau : (W, w_0) \to (W', w'_0)$ 
is called a {\it Legendre diffeomorphism-germ} 
if $\tau$ maps $\pi$-fibers to $\pi'$-fibers, or more 
exactly, if there exists a diffeomorphism-germ 
$\bar{\tau} : (Z, z_0) \to (Z', z'_0)$ 
such that $\bar{\tau}\circ \pi = \pi'\circ\tau$. 

A pair $(f, \pi)$ is {\it Legendre equivalent} to 
$(f', \pi')$ if there exists a contact equivalence $(\sigma, \tau)$ 
of $f$ and $f'$ such that $\tau$ is a Legendre diffeomorphism. 
In this case, we call $(\sigma, \tau)$ a {\it Legendre 
equivalence} of $(f, \pi)$ and $(f', \pi')$. 

An integral map-germ $f : (N, x_0) \to W$ 
is called {\it homotopically Legendre stable} 
if any integral deformation $(f_t)$ of $f$ is 
trivialized under Legendre equivalence: 
$$
\tau_t\circ f_t\circ \sigma_t^{-1} = f, 
$$
$(\sigma_t, \tau_t)$ being Legendre equivalences of $f_t$ and $f$. 
Here $\sigma_t$ may move base points of germs.

Moreover we can define, over the $\R$, 
the notion of Legendre stability of 
map-germs: 
Roughly speaking, 
an integral map-germ $f : (N, x_0) 
\to W$ is Legendre stable with 
respect to an Legendre fibration 
$\pi : W \to Z$ if, by any sufficiently small 
integral perturbations, the Legendre equivalence class of $(f, \pi)$ 
is not removed. 
To formulate accurately, 
denote by $C^{\infty}_I(N, W)$ the space of $C^\infty$ integral mappings 
from $N$ to $W$, endowed with the Whitney $C^\infty$ topology. 
Then an integral map-germ $f : (N, x_0) \to W$ is {\it Legendre   
stable} if, 
for any integral representative $f : U \to W$ of $f$, there exists 
a neighborhood $\Omega$ in 
$C^{\infty}_I(U, W)$ of $f$ 
such that, for any $f' \in \Omega$, the 
original pair $(f_{x_0}, \pi)$ of germs 
is Legendre equivalent to $(f'_{x'_0}, \pi)$ for some 
$x'_0 \in U$ (cf. \cite{AGV}). 

To characterize the Legendre stability by means of 
transversality, we introduce the notion of  
{\it integral jet spaces}. 
Denote by $J^r_I(N, W)$ the set of $r$-jets 
of integral map-germs $f : (N, x_0) \to (W, w_0)$ {\it 
of corank at most one}:  
$$
J^r_I(N, W) = \{ j^rf(x_0) \mid f : (N, x_0) \to (W, w_0) 
{\mbox {\rm \ integral, corank}}_{x_0}f \leq 1\}. 
$$
Then $J^r_I(N, W)$ is 
a submanifold of the ordinary jet space $J^r(N, W)$ 
(\S \ref{Integral jet}). 
Moreover, for $j^rf(x_0) \in J^r_I(N, W)$, 
the {\it Legendre equivalence class} of $j^rf(x_0)$, namely, 
the set of 
$r$-jets of map-germs which are Legendre equivalent to $f : 
(N, x_0) \to (W, w_0)$ form a submanifold of $J^r_I(N, W)$. 

If $f : N \to W$ is an integral mapping of corank at most one, then 
the image of the $r$-jet extension $j^rf : N \to J^r(N, W)$ 
is contained in $J^r_I(N, W)$. Then we regard $j^rf$ as a mapping to 
$J^r_I(N, W)$. 
Based on a Legendre version of 
transversality theorem (\S \ref{Integral jet}), 
Legendre stability is characterized by the transversality. 

\smallskip 

We apply, over $\C$, the transversality as the definition of stability. 

For a manifold-germ $(N, x_0)$, we denote by 
${\mathcal E}_{N, x_0}$ the $\KK$-algebra consisting of $C^\infty$ function-germs 
$(N, x_0) \to \KK$, and by $m_{N, x_0}$ the unique maximal ideal of  
${\mathcal E}_{N, x_0}$. If the base point $x_0$ is clear in the context, we abbreviate 
${\mathcal E}_{N, x_0}$ and $m_{N, x_0}$ to ${\mathcal E}_N$ and $m_N$ 
respectively. 

Now set 
$$
r_0 = {\mbox{\rm inf}}\{ r \in \NN \mid f^*{\mathcal E}_{W}\cap m_N^{r+1} 
\subset f^*m_{W}^{n+2}\}. 
$$
If 
$f : (N, x_0) \to W$ is an open Whitney umbrella, then $f$ is, in particular, 
finite, namely, ${\mathcal E}_N$ is a finite ${\mathcal E}_W$-module via $f^* : {\mathcal E}_W \to {\mathcal E}_N$. 
Therefore
 $r_0$ is a finite positive integer, determined by $n$ and $k$, the type of the 
open Whitney umbrella. 
Actually $r_0$ depends only on the right-left equivalent class of $f$. 

The main purpose of the present paper is 
to show the following: 

\bet
\label{Stability} 
{\mbox{\rm (}}Arnol'd-Mather type 
characterization of Legendre stability{\mbox{\rm )}}. 
For an integral map-germ 
$f : (N^n, x_0) 
\to (W^{2n+1}, w_0)$ of corank at most one, 
the following conditions are equivalent to each other: 

{\mbox{\rm (s)}}
 $f$ is Legendre stable. 

{\mbox{\rm (hs)}}
 $f$ is homotopically Legendre stable. 

{\mbox{\rm (is)}}
 $f$ is infinitesimally Legendre stable. 

{\mbox {\rm (a)}} 
$f$ is an open Whitney umbrella and 
$f^*{\mathcal E}_{W}$ is generated by $1, p_1\circ f, \dots, p_n\circ f$ 
as ${\mathcal E}_Z$-module via $(\pi\circ f)^*$.

{\mbox {\rm (a$'$)}} 
 $f$ is an open Whitney umbrella and 
 $Q(f) := f^*{\mathcal E}_W/(\pi\circ f)^*m_Z{\mathcal E}_W$ 
is generated over $\KK$ by 
 $1, p_1\circ f, \dots, p_n\circ f$. 
 
{\mbox {\rm (a$''_r$)}} $(r \geq r_0)$. 
 $f$ is an open Whitney umbrella and 
$Q_{r+1}(f) := f^*{\mathcal E}_{W}/\{(\pi\circ f)^*m_Zf^*{\mathcal E}_{W} + 
f^*{\mathcal E}_{W} \cap m^{r+2}_N \}$ is 
generated by $1, p_1\circ f, \dots, p_n\circ f$ over $\KK$. 

{\mbox {\rm (t$_r$)}} $(r \geq r_0)$.  The jet extension $j^rf : (N, x_0) \to 
J_I^r(N, W)$ is transversal to the Lagrange 
equivalence class of $j^rf(x_0)$. 

\ent

We must explain the notion of infinitesimal Legendre stablity (is): 
Of course, it is the infinitesimal counterpart of Legendre stability. 
Now recall the notion of infinitesimal stability due to Mather \cite{Mather} 
for a general $C^\infty$ map-germ $f : (N, x_0) \to W$. 
A map-germ $f$ is called infinitesimally stable if 
$V_f = tf(V_N) + wf(V_W)$, where 
$V_N$ (resp. $V_W$, $V_f$) 
is the module consisting of all germs of vector fields over $(N, x_0)$ 
(resp. over $(W, f(x_0))$, along $f$), and 
$tf : V_N \to V_f$ (resp. $wf : V_W \to V_f$) 
is defined by $tf(\xi)(x) = f_*(\xi(x)), (\xi \in V_N, x \in (N, x_0))$ 
(resp. $wf(\eta)(x) = \eta(f(x)), (\eta \in V_W, x \in (N, x_0))$). 
Similarly we call an integral map-germ $f : (N, x_0) \to W$ 
{\it infinitesimally Legendre stable} if 
$VI_f = tf(V_N) + wf(VL_W)$, 
where $VL_W (\subset V_W)$ (resp. $VI_f (\subset V_f)$) is the module 
of all germs of infinitesimal Legendre deformations over $(N, x_0)$ 
(resp. infinitesimal integral deformations of $f$). 
See \S \ref{Infinitesimal deformations.}. 

\ 

The equivalence of {\mbox{\rm (hs)}}
and {\mbox{\rm (is)}} is one of consequences of 
{\it Legendre versality theorem}: 
We introduce the notion of 
Legendre versality of integral deformations of 
integral map-germs. 

A deformation $F : (N\times \KK^r, (x_0, 0)) 
\to W$ of an integral map-germ 
$f : (N, x_0) \to W$ is called {\it integral} 
if each $f_\lambda = F\vert_{N\times\{\lambda\}}, (\lambda \in (\KK^r, 0))$ 
is integral, for a representative of $F$. We write $F = (f_\lambda)$ in short. 
An integral deformation $F$ of $f$ is called {\it Legendre versal} if 
any other integral deformation $G : (N\times \KK^s, (x_0, 0)) 
\to W$ of $f$ is induced from $F$ up to Legendre equivalence, 
namely if there exist a map-germ $\varphi : (\KK^s, 0) \to (\KK^r, 0)$ 
and a Legendre deformation $(\sigma_\mu, \tau_\mu), (\mu \in (\KK^s, 0))$ 
such that $g_\mu = \tau_\mu\circ f_{\varphi(\mu)}\circ\sigma_\mu^{-1}$ 
for any $(\mu \in (\KK^s, 0))$, where $g_\mu(x) = G(x, \mu)$. 
$F$ is called {\it infinitesimally Legendre versal} 
if 
$$
VI_f = \left.\left.\left\langle \dfrac{\pa F}{\pa \lambda_1}\right\vert_{\lambda = 0}, 
\dots, \dfrac{\pa F}{\pa \lambda_r}\right\vert_{\lambda = 0}
\right\rangle_{\KK} + tf(V_N) + wf(VL_W). 
$$

Then we also mention in this paper on a proof of the following: 

\bet
\label{Versality}
An integral deformation $F : (N\times 
\KK^r, (x_0,0)) \to W$ of an integral 
map-germ $f : (N, x_0) \to W$ 
of corank at most one,  
is Legendre versal 
if and only if $F$ is infinitesimally Legendre versal. 
Any Legendre versal deformations of $f$ with 
the same number of parameters are Legendre 
equivalent to each other. 
An integral map-germ $f : (N, x_0) \to W$ 
of corank at most one
has a Legendre versal deformation if and only if 
$tf(V_N) + wf(VL_W)$ is of finite codimension over $\KK$ in $VI_f$. 
\ent

Setting $r = 0$ we have again that 
$f$ is homotopically Legendre stable if and only if 
$f$ is infinitesimally Legendre stable.

\section{Contact stability}
\label{contact stability}

Related to the notion of Legendre stability, 
we define the notion of contact stability of 
map-germs: 
Roughly speaking, 
an integral map-germ $f : (X, x_0) 
\to W$ is {\it contact stable} if, by any sufficiently small 
integral perturbations, the contact equivalence class of $f_{x_0}$ 
is not removed but remains nearby $x_0$. 

More exactly, an integral map-germ $f : (N, x_0) \to W$ is {\it contact  
stable} if, 
for any integral representative $f : U \to W$ of $f$, there exists 
a neighborhood $\Omega$ in 
$C^{\infty}_I(N, W)$ such that, for any $f' \in \Omega$, the 
original germ $f$ 
is contact equivalent to $f'_{x'_0}$ for some 
$x'_0 \in U$ (cf. \cite{AGV} page 325). 

An integral map-germ 
$f : (N, x_0) \to W$ is called {\it homotopically contact stable} 
if any one-parameter integral deformation $F = (f_t)$ of $f$ is trivialized 
by contactomorphisms: 
$$
\tau_t\circ f_t\circ \sigma_t^{-1} = f, 
$$
$(\sigma_t, \tau_t)$ being contactomorphism of $f_t$ and $f$. 
Here $\sigma_t$ may move base points of germs. 
(See the definition of contactomorphism in \S 1). 
$f$ is called 
{\it infinitesimally contact stable} if $f$ satisfies the 
infinitesimal condition corresponding to the contact stability, 
namely, if $f$ satisfies the condition 
$$
VI_f = tf(V_N) + wf(VH_W). 
$$
For a map-germ $f : (N, x_0) \to (W, w_0)$ we set 
$$
R_f := \{ h \in {\mathcal E}_W 
\mid dh \in {\mathcal E}_Nd(f^*{\mathcal E}_W) \}, 
$$
for the exterior differential $d$. Here we denote by ${\mathcal E}_W$ and 
${\mathcal E}_N$ the algebra of function-germs on $W$ and $N$ respectively.

Then we have: 

\bep
{\mbox{\rm (}}Classification of 
contact stable germs{\mbox{\rm )}}. 
\label{contact stable}
Let 
$f : (N, x_0) \to W^{2n+1}$ be 
an integral map-germs of corank at most one. 
Then the following conditions are equivalent: 

{\mbox{\rm (cs)}}
$f$ is contact stable.  

{\mbox{\rm (hcs)}}
$f$ is homotopically contact stable. 

{\mbox{\rm (ics)}}
$f$ is infinitesimally contact stable.

{\mbox{\rm (owu)}}
$f$ is an open Whitney umbrella. 

{\mbox{\rm (ca)}}
$R_f = f^*{\mathcal E}_W$ and 
$f$ is diffeomorphic {\rm (}i.e. ${\mathcal A}$-equivalent{\rm )} 
to an analytic map-germ $f' : (\KK^n, 0) 
\to (\KK^{2n+1}, 0)$ {\rm (}not necessarily integral{\rm )} 
such that the codimension of the 
singular locus of the complexification $f'_{\C}$ of $f'$ is greater than 
or equal to $2$.

{\mbox{\rm (ct)}}
The jet extension $j^rf : (N, x_0) \to J^r_I(N, W)$ is 
transversal to the contactomorphism class of $j^rf(x_0)$, 
for an integer $r \geq \dfrac{n}{2} + 1$. 

\enp

The notions in Proposition \ref{contact stable} 
are discussed in detail along the following sections, 
in particular in \S \ref{Infinitesimal deformations.} and 
in \S \ref{Integral jet}. 
The proof of Proposition \ref{contact stable} 
will be given in \S \ref{Legendre}. 

Similarly to Legendre versality theorem (Theorem \ref{Versality}), 
we can show {\it contact versality theorem}, which gives an alternative 
proof of the equivalence of (hcs) and (ics). 

An integral deformation $F = (f_\lambda) : 
(N\times \KK^r, (x_0, 0)) 
\to W$ of an integral map-germ $f : (N, x_0) 
\to W$ is called {\it contact versal} if 
any other integral deformation $G = (g_\mu) : (N\times \KK^s, (x_0, 0)) 
\to W$ of $f$ is induced from $F$ up to contactomorphisms, 
namely if there exist a map-germ $\varphi : (\KK^s, 0) \to (\KK^r, 0)$ 
and a family of contactomorphisms  
$(\sigma_\mu, \tau_\mu), (\mu \in (\KK^s, 0))$ 
such that $g_\mu = \tau_\mu\circ f_{\varphi(\mu)}\circ\sigma_\mu^{-1}$ 
for any $(\mu \in (\KK^s, 0))$, where $g_\mu(x) = G(x, \mu)$.

$F$ is called {\it infinitesimally contact versal} if 
$$
VI_f = \left\langle \left.\dfrac{\pa F}{\pa \lambda_1}\right\vert_{\lambda = 0}, 
\dots, \left.\dfrac{\pa F}{\pa \lambda_r}\right\vert_{\lambda = 0}
\right\rangle_{\R} + tf(V_N) + wf(VH_W). 
$$

\bet
\label{contact Versality}
An integral deformation $F : (N\times 
\KK^r, (x_0,0)) \to W$ of an integral 
map-germ $f : (N, x_0) \to W$ 
of corank at most one,  
is contact versal 
if and only if $F$ is infinitesimally contact versal. 
Any contact versal deformations of $f$ with 
the same number of parameters are contactomorphic  
to each other. 
An integral map-germ $f : (N, x_0) \to W$ 
of corank at most one
has a contact versal deformation if and only if 
$tf(V_N) + wf(VH_W)$ is of finite codimension over $\KK$ in $VI_f$. 
\ent

\section{Lie derivative.} 
\label{Lie derivative.}

Let $N, W$ be manifolds, and $f : N \to W$ a mapping. 
A mapping $v : N \to TW$ is called a {\it vector field along} $f$ or 
an {\it infinitesimal deformation} of $f$, 
if $\pi\circ v = f$, for 
the projection $\pi : TW \to W$. 
We denote by $V_f$ the module of all vector field along $f$. 
By the fiberwise addition and scalar multiplication on $TW$, 
$V_f$ turns out to be a module over the function-algebra ${\mathcal E}_N$ 
on $N$. 

It is easy to see that there exists a one-parameter 
deformation $F : U \to W$ of $f$ defined on an open neighborhood $U$ 
in $N \times \KK$ of $N \times \{ 0\} 
\cong N$ such that $F\vert_{N\times\{ 0\}} = f$. 
We write as $F = (F_t)$ so that $F_0 = f$. 
Then we define, for a differential $p$-form $\alpha$ on $W$, 
a differential $p$-form $L_v\alpha$ on $N$ by 
$$
L_v\alpha = \left.\dfrac{d}{dt}\right\vert_{t=0} F_t^*\alpha. 
$$
For this, see also \cite{Inv} p.225. 
Then 
$L_v\alpha$ does not depend on the choice of $F$ but depends only on $v$. 
We call $L_v\alpha$ the {\it Lie derivative} of $\alpha$ by $v$. 
Moreover we define the interior product $i_v\alpha$, that is 
a differential $(p-1)$-form on $N$ by 
$$
i_v\alpha(Z_1, \dots, Z_{p-1})(x) = \alpha(v(x), 
f_*Z_1(x), \dots, f_*Z_{p-1}(x)), 
$$
for vector fields $Z_1, \dots, Z_{p-1}$ over $N$. 

\smallskip 

\bee
{\rm
Let $N = TW$ and $f = \pi : TW \to W$. We 
regard the identity map $1 : TW \to TW$ as 
a vector field along $\pi$. Then, for a $p$-form $\alpha$ 
on $W$, we have defined 
the $p$-form $L_1\alpha$ and $(p-1)$-form $i_1\alpha$ on $TW$. 
}\QED
\ene

\bel
\label{i and L}
We have the following fundamental formulae: 
$$
\begin{array}{lrcl}
{\mbox{\rm (1)}} & 
i_v(\lambda\alpha + \mu\beta) 
& = & \lambda(i_v\alpha) + \mu(i_v\beta), \\
{\mbox{\rm (2)}} & 
i_{\lambda u + \mu v}\alpha 
& = & (f^*\lambda)(i_u\alpha) + (f^*\mu)(i_v\alpha), \\
{\mbox{\rm (3)}} &  
L_v\alpha & = & i_v(d\alpha) + d(i_v\alpha),  \\
{\mbox{\rm (4)}} & 
L_v(\alpha \wedge \beta) & = & (L_v\alpha) \wedge f^*\beta 
+ f^*\alpha \wedge (L_v\beta),  \\
{\mbox{\rm (5)}} & 
i_v(\alpha \wedge \beta) & = & (i_v\alpha) \wedge f^*\beta + 
(-1)^r f^*\alpha \wedge (i_v\beta). 
\end{array}
$$
Here $u, v$ are vector fields along a mapping $f : N \to W$, 
$\lambda, \mu$ are functions on $W$, 
$\alpha, \beta$ are differential forms on $W$, 
and $\alpha$ is an $r$-form. 

In particular, we refer (3) as 
the Cartan type formula: $L_v = di_v + i_vd$. 
\enl

\Proof 
(1) and (2) are straightforward from the definition. 
The proof of (3) 
is given in \cite{Inv} Lemma 3.3. (4), (5) are easily proved similarly to 
the ordinary case $W = N$ and $f$ is the identity mapping. 
\QED

The following formulae are proved from the definitions in the straightforward way. 

\bel
\label{pullback.etc}
Let $f : N \to W$ be a mapping, 
$v : N' \to TN$ a vector field along a mapping $N' \to N$, 
$w : W \to TW'$ a vector field along a mapping $W \to W'$, 
$\alpha$ a differential form on $W$ and 
$\alpha'$ a differential form on $W'$. Then we have 
$$
\begin{array}{lcc}
{\mbox{\rm (i)}} & 
L_{w\circ f}\alpha' = f^*(L_w\alpha'),  & \quad 
L_{f_*v}\alpha = L_v(f^*\alpha),  \\
{\mbox{\rm (ii)}} & 
i_{w\circ f}\alpha' = f^*(i_Y\alpha'), & \quad 
i_{f_*v}\alpha = i_v(f^*\alpha). 
\end{array}
$$
Here $w\circ f$ is the pull-back of $w$ by $f$, and 
$f_*v$ is the push-forward of $v$ by $f$: $(w\circ f)(x) = w(f(x)), (x \in N)$, 
$(f_*v)(x') = f_*(v(x')), (x' \in N')$. 

In particular we have  
$$
\begin{array}{lcc}
{\mbox{\rm (i')}} & 
L_{Y\circ f}\alpha = f^*(L_Y\alpha),  & \quad 
L_{f_*X}\alpha = L_X(f^*\alpha),  \\
{\mbox{\rm (ii')}} & 
i_{Y\circ f}\alpha = f^*(i_Y\alpha), & \quad 
i_{f_*X}\alpha = i_X(f^*\alpha), 
\end{array}
$$
for a vector field $X$ over $N$, 
 a vector field $Y$ over $W$, and for a differential form $\alpha$ on $W$. 
Here $Y\circ f$ is the pull-back of $Y$ by $f$, and 
$f_*X$ is the push-forward of $X$ by $f$. 
\enl

\ 

Then the fundamental concept of this paper is introduced as follows: 
\bep 
\label{natural lifting}
Let $W$ be a manifold and $\alpha$ a differential form on $W$. 

{\mbox{\rm (1)}} 
There exists a unique differential form $\widetilde{\alpha}$ on 
$TW$ such that, for any vector field $X : W \to TW$ 
over $W$, $X^*{\widetilde{\alpha}} = {L}_X\alpha$ holds. 

{\mbox{\rm (2)}} 
Moreover, $\widetilde{\alpha}$ of {\mbox{\rm (1)}} satisfies 
$v^*{\widetilde{\alpha}} = {L}_v\alpha$, for any vector field 
$v : N \to TW$ along a mapping $f : N \to W$. 

{\mbox{\rm (3)}} 
$d\widetilde{\alpha} = \widetilde{d\alpha}$ and 
$\widetilde{f^*\alpha} = (f_*)^*(\widetilde\alpha)$, 
where $f_* : TN \to TM$ 
is the bundle homomorphism 
defined by differential of $f$. 

\enp 

In fact, we have $\widetilde{\alpha} = L_1\alpha$. 
We call $\widetilde{\alpha}$ the {\it natural lifting} of $\alpha$. 
The notion of natural liftings is first defined, even for general tensors, in 
\cite{YI}\cite{YK} in a different manner: This fact is pointed out to 
the author by H. Sato. Though our construction is limited to differential 
forms, it seems more direct and useful for the infinitesimal study of 
differential systems. We are going to apply, in this paper, 
 the notion of natural liftings for the infinitesimal study of stability of 
integral mappings in contact geometry. 

\smallskip

\noindent{\it Proof of Proposition \ref{natural lifting}:} 
(1) 
We set $\widetilde{\alpha} = L_1\alpha$, for the identity mapping 
$1 : TW \to TW$. Then $X^*\widetilde{\alpha} = X^*L_1\alpha = 
L_{1\circ X}\alpha = L_X\alpha$. 
Similarly we have (2). 
Let, for another $\beta$, $X^*\beta = X^*\widetilde{\alpha}$, for any vector 
field over $W$. Then, for any $z \in TW$ and any $v \in T_z(TW) \setminus K$, 
there exists a vector field $X$ over $W$ and $u \in T_{\pi(z)}W$ 
such that $X_*(u) = v$. Here $\pi : TW \to W$ the canonical 
projection and $K$ is the kernel of $\pi_* : T_z(TW) \to T_{\pi(z)}W$. 
Then $\langle \beta, v\rangle = \langle X^*\beta, u\rangle = 
\langle X^*\widetilde{\alpha}, u\rangle = \langle \widetilde{\alpha}, v\rangle$. 
Thus $\beta$ and $\widetilde{\alpha}$ coincide on $T_z(TW) \setminus K$ thus on 
$T_z(TW)$, the linear-hull of $T_z(TW) \setminus K$, 
for any $z \in TW$. Therefore $\beta = \widetilde{\alpha}$. 
(3) follows from the uniqueness of the natural lifting of $d\alpha$ and 
$f^*\alpha$: For example, 
$X^*(f_*)^*(\widetilde{\alpha}) = (f_*X)^*(\widetilde{\alpha}) = 
L_{f_*X}\alpha = L_X(f^*\alpha)$, for any vector field $X$ over $N$.  
\QED

\smallskip

\bee
{\rm
Let $M$ a symplectic manifold, and $\omega$ 
the symplectic form on $M$. Since $\omega$ is 
non-degenerate, $\omega$ induces an 
isomorphism $TM \cong T^*M$. On the other hand, 
$T^*M$ is endowed with the canonical 
symplectic form $d\theta_M$, which 
is independent of the symplectic structure of $M$. 
Therefore $d\theta_M$ is regarded as 
a symplectic form on $TM$. This coincides with the 
natural lifting 
$\widetilde{\omega}$.  
}
\ene

\smallskip  

\bee
{\rm
Let $(p, q, r)$ be a 
Darboux coordinates of $(W, D)$ at a point $w_0 \in W$. 
Then the standard contact 
form $\alpha = dr - pdq$ gives the contact distribution $D \subset TW$. 
Let $(p, q, r; \phi, \xi, s)$ be the induced local coordinates of the 
tangent bundle $TW$; $(\phi, \xi, s)$ being fiber coordinates. 
Then we have 
$$
\widetilde{\alpha} = d(s - p\xi) + \xi dp - \phi dq.
$$
Remark that $\widetilde{\alpha}$ is linear in the 
fiber coordinates $(\phi, \xi, s)$. 
}
\ene

\smallskip 

In general we have 
\bel
\label{additive} 
Let $f : N \to W$ be a mapping, and $\alpha$ a differential form on $W$. 
Then, 
for $v_1, v_2 \in V_f$, we have 
$i_{v_1 + v_2}\alpha = i_{v_1}\alpha + i_{v_2}\alpha$, and 
$(v_1 + v_2)^*\widetilde{\alpha} = 
v_1^*\widetilde{\alpha} + v_2^*\widetilde{\alpha}$. 
\enl

\Proof
The first equality follows from the definition of interior product. 
The second equality follows from  
Proposition \ref{natural lifting} (2) and the Cartan's formula $L_v = di_v + i_vd$. 
\QED

\ 

The notion of natural liftings is defined also for differential systems. 
Let $W$ be a manifold and $\Omega$ the sheaf of differential 
forms on $W$. A subsheaf $I \subset \Omega$ 
is called a {\it differential system} on $W$ 
if it is a $d$-closed ideal of the differential 
algebra $\Omega$, namely, if, 
for any section $\alpha$ of $I$ and for any section $\beta$ of $\Omega$ 
(defined on the same open subset of $W$), 
$\alpha\wedge\beta$ and $d\alpha$ are sections of $I$. 
Let $S$ be a set 
of differential forms on open subsets of $W$. Then 
the differential system $\langle S\rangle$ 
generated by $S$ has the stalk $\langle S\rangle_x$, for each $x \in W$, 
consisting of the functional linear combination of elements $\alpha_x\wedge\beta_x$  
and $d\alpha_x\wedge\beta_x$, for those $\alpha \in S$ and 
differential forms $\beta$ defined over $x$. 

For example, a contact structure $D \subset TW$ on $W$ may be defined also as 
the differential system generated by local sections of $D^{\perp} \subset 
T^*W$, local contact forms 
compatible with $D$. 

Let $I$ be a differential system on $W$. Then 
the {\it natural lifting} $\widetilde{I}$ of $I$ is defined as 
the differential system on $TW$ generated by 
the natural liftings $\widetilde{\alpha}$ of all sections $\alpha$ 
of $I$. 
If $f : N \to W$ is a mapping, then $f^*I$ denotes the differential system 
generated by $f^*\alpha$ for all sections of $I$.  
Then we have by Proposition \ref{natural lifting} (3): 
\bel
Let $I$ be a differential system on $W$.   Then $\widetilde{f^*I} = (f_*)^*(\widetilde{I})$, 
where $f_* : TN \to TW$ is the differential mapping of $f$. 
\enl

\section{Contact Hamilton vector fields.}
\label{Contact Hamilton vector fields.}

Let $(W, D)$ be a contact manifold, and 
$\alpha$ a local contact form representing $D$. 
There does not necessarily exist $\alpha$ globally; $\alpha$ can be taken 
over an open subset of $W$ 
where the contact distribution $D$ is co-oriented. 
A vector field $X$ over $W$ is called a {\it contact 
vector field} if the Lie derivative $L_X\alpha = \mu\alpha$ 
for a function $\mu$, namely if 
$X$ preserves the contact distribution $D$. 

Deleting $W$ if necessary, 
we assume a contact form $\alpha$ is taken over $W$. 
Let $H : W \to \KK$ be a function. 
Then there exists a unique contact vector field $X = X_H$ over $W$ 
with the condition $i_X\alpha = H$. 
The contact vector field $X_H$ is called the {\it contact 
Hamilton vector field with Hamilton function} $H$. 

If $\alpha = dr - \sum_{i=1}^n p_idq_i$, then $X_H$ is 
explicitly given by 
$$
{X_H} = 
\sum_{i=1}^n \left( \dfrac{\pa H}{\pa q_i} + 
p_i\dfrac{\pa H}{\pa r}\right)\dfrac{\pa}{\pa p_i} 
- \sum_{i=1}^n \dfrac{\pa H}{\pa p_i}\dfrac{\pa}{\pa q_i} + 
\left( H - \sum_{i=1}^n p_i\dfrac{\pa H}{\pa p_i}\right)\dfrac{\pa}{\pa r}. 
$$

Conversely, any contact vector field is locally a contact Hamilton 
vector field with some Hamiltonian function. 

Associated to a contact form $\alpha$, we define the {\it Reeb vector 
field} $R$ by $i_R\alpha = 1, i_Rd\alpha = 0$. 
Note that, since $\alpha$ is a contact form, $R$ is characterised uniquely. 
If $\alpha = dr - pdq$, then $R = \dfrac{\pa}{\pa r}$. 
Then we have: 
\bel
\label{LieReeb}
Let $\alpha$ be a contact form on $W$, and $H : W \to \KK$ a 
function. Then we have 

{\rm (1)} 
$L_{X_H}\alpha = R(H)\alpha$ and 
$i_{X_H}d\alpha = R(H)\alpha - dH$. 

{\rm (2)} 
Let $\eta$ be a vector field on $W$. If 
$i_{\eta}d\alpha = 0$, then $\eta = (i_{\eta}\alpha)R$. 

{\rm (3)} 
$X_1 = R$. 
\enl

\Proof
(1) The first equality holds for a system of 
coordinates $(p, q, r)$ with $\alpha = dr - pdq$. 
Remark that $X_H$ and $R$ are defined intrinsically from the contact 
form $\alpha$. 
The latter equality follows from 
$L_{X_H}\alpha = i_{X_H}d\alpha + di_{X_H}\alpha = 
i_{X_H}d\alpha + dH$. 
(2) Set $\eta' = \eta - (i_{\eta}\alpha)R$. Then 
$i_{\eta'}d\alpha = 0$ and $i_{\eta'}\alpha = 0$. 
Therefore $\eta' = 0$, and we have $\eta = (i_{\eta}\alpha)R$. 
(3) By (1), we have $i_{X_1} d\alpha = 0$. Since $i_{X_1}\alpha = 1$, we see 
$X_1 = R$. 
\QED

\smallskip

We have the following formula for 
the contact Hamilton vector field with 
the sum (resp. product) of two contact Hamilton functions: 

\bel
\label{KH}
For functions $K, H$ on $W$, we have 
$$ 
X_{K + H} = X_K + X_H, 
$$
$$
X_{KH} = K\cdot X_H + H\cdot X_K - (KH)\cdot R = 
K\cdot X_H + H\cdot X_K - (KH)\cdot X_1. 
$$
In particular, $X_{aH} = aX_H, (a \in \KK)$. 
\enl

\Proof
The first one is clear. To show the second equality, we 
set $\eta = K\cdot X_H + H\cdot X_K - X_{KH}$. Then 
$$
\begin{array}{rcl}
i_{\eta}d\alpha & = & K(R(H)\alpha - dH) + H(R(K)\alpha - dK) 
- (R(KH)\alpha - d(KH)) \\
 & = & (KR(H) + HR(K) - R(KH))\alpha = 0. 
\end{array} 
$$
Moreover, 
$i_{\eta}\alpha = KH + HK - KH = KH$. 
Therefore, by Lemma \ref{LieReeb}, 
$\eta = (KH)\cdot R = (KH)\cdot X_1$. 
\QED

\ 

We denote by $VH_W$ the vector space of 
contact Hamilton vector fields over $W$ and by ${\mathcal E}_W$ 
the $\KK$-algebra of functions on $W$. 
Define a linear map $\Phi : {\mathcal E}_W \to VH_W$ by 
$\Phi(H) = X_H$. Then $\Phi$ 
is an isomorphism of vector spaces. Therefore $VH_W$ 
is endowed with ${\mathcal E}_W$-module structure 
induced from $\Phi$, namely, $K*X_H = X_{KH}$. 
Here, we distinguish this new functional multiplication, 
using $*$, with the ordinary functional multiplication in 
$V_W$, the ${\mathcal E}_W$-module consisting of all vector fields over $W$. 

\

In term of the local coordinates $p, q, r$ of $(W, w_0)$ with $\alpha = 
dr - \sum_{i=1}^n p_idq_i$, 
we define the order of function-germs $h = h(p, q, r) \in {\mathcal E}_W$ 
by setting 
$$
\weight(p_i) = \weight(q_j) = 1, (1\leq i, j \leq n), \quad 
{\mbox{\rm and\ }} \weight(r) = 2; 
$$
namely, 
$\ord(h) \geq r$ if the Taylor expansion of $h$ has no monomials of weight $< r$. 
We set $m_W^{(r)} := \{ h \in {\mathcal E}_W \mid 
\ord(h) \geq r \}$. 

If $\tau : (W, w_0) \to (W, w_0)$ is a contactomorphism, then 
$\ord(h\circ \tau) = \ord(h)$. 
Then we can define, on the local ring ${\mathcal E}_W$, the filtration 
$$
{\mathcal E}_W \supset m_W^{(1)} \supset  m_W^{(2)} \supset \cdots 
\supset  m_W^{(r)} \supset \cdots. 
$$
Note that 
$$
m_W^{(2r)} \subset m_W^r \subseteq m_W^{(r)}, (r = 0, 1, 2, \dots). 
$$
In particular 
$m_W^2 \subset m_W^{(2)} \subset m_W$. 

In the ${\mathcal E}_W$-module $VH_W$ introduced above, 
we have 
$$ 
m_W^2*VH_W \subseteq VH_W \cap m_WV_W = m_W^{(2)}*VH_W. 
$$

\

Let $\pi : W \to Z$ be a Legendre fibration. Then 
a contact vector field $X$ over $W$ is called a {\it 
Legendre vector field} if, $X$ is lowerable, namely, 
if there exists a vector field 
$Y$ over $Z$ such that $t\pi(X) = w\pi(Y)$ as vector fields along $\pi$. 
Then easily we have: 
\bep
\label{affine function}
Let $(p_1, \dots, p_n, q_1, \dots, q_n, r)$ 
be a Darboux coordinate, so that 
$\alpha = dr - pdq$. Then 
a contact Hamilton vector field $X_H$ with 
Hamilton $H = H(p, q, r)$ is a Legendre 
vector field if and only if 
$H$ is an affine function, namely, $H$ 
is of form 
$$
H(p, q, r) = a_0(q, r) + a_1(q, r)p_1 + \cdots + 
a_n(q, r)p_n
$$
\enp

We denote by $VL_W = VL_{(W, \pi)}$,  
the totality of Legendre vector fields over $W$ with 
respect to $\pi$.

\section{Infinitesimal deformations.} 
\label{Infinitesimal deformations.} 

Let $f : (N, x_0) \to W$ be an integral map-germ. 
The space of infinitesimal integral deformations 
of $f$ is, at least formally, given by 
$$
VI_f = 
\{ v : (N, x_0) \to TW 
\mid v^*\widetilde{\alpha} = 0, \ \pi\circ v = f \}, 
$$
where $\pi : TW \to W$ is the natural projection, and 
$\widetilde{\alpha}$ is the natural lifting to $TW$ of a contact 
$1$-form $\alpha$ locally defining $D$ near $w_0 = f(x_0) \in W$. 

Recall that $VH_W$ denotes the ${\mathcal E}_W$-module 
of contact Hamilton vector fields over $W$. 
Define a linear mapping 
$wf : VH_W \to VI_f$ by $wf(H) = X_H\circ f, (H \in {\mathcal E}_W)$. 

For $v \in VI_f$, we call $i_v\alpha \in {\mathcal E}_N$ the {\it 
generating function} of $v$. The linear mapping 
$e : VI_f \to R_f$ is defined by taking generating function. 
Here 
$$
R_f := \{ h \in {\mathcal E}_N \mid dh \in {\mathcal E}_Nd(f^*{\mathcal E}_W) \}. 
$$

In local coordinates, we have $e(v) = s\circ v - \sum (p\circ f)(\xi\circ v)$ and 
$$
0 = v^*\widetilde{\alpha} = 
d(e(v)) + \sum (\xi\circ v)d(p\circ f) - 
\sum (\phi\circ v)d(q\circ f). 
$$
Therefore $e(v) \in R_f$. 
Note that 
$i_v(\lambda\alpha) = (\lambda\circ f) i_v\alpha$. 

We see the mapping $e$ is surjective. In fact, for any $h \in R_f$, $dh$ is a 
functional linear combination of the exterior derivatives of 
components of $f$. Since $f$ is integral, $r\circ f$ is a 
functional linear combination of $d(p\circ f), d(q\circ f)$, and so is $dh$. 
Therefore, choosing $\xi\circ v, \varphi\circ v$ and $s\circ v$ properly, 
we get $v \in VI_f$ with $e(v) = h$. 

Note that 

\bel
\label{generating function 1}
We have $i_{X_H\circ f}\alpha = f^*(i_{X_H}\alpha) = f^*H$. 
Therefore the generating function of $X_H\circ f$ is equal to the 
pull-back $f^*H$ of the Hamiltonian function $H$. 
\enl

We need a result proved in page 222 of \cite{Inv}: 

\bel
\label{R_f finite}
Let $f : (N, x_0) \to W$ be of corank $\leq 1$. 
If 
$$
R_f := \{ e \in {\mathcal E}_N 
\mid de \in {\mathcal E}_Nd(f^*{\mathcal E}_W) \}
$$ is a finite ${\mathcal E}_W$-module if and only if $f$ is a finite map-germ, 
namely, ${\mathcal E}_N$ is a finite ${\mathcal E}_W$-module via $f^* 
: {\mathcal E}_W \to {\mathcal E}_N$. 
\enl

Now set 
$VI_f' = {\mbox{\rm Ker}}(e : VI_f \to R_f)$. 
Then we have the exact sequence of vector spaces: 
$$
0 \lon VI_f' \lon VI_f \stackrel{e}{\lon} R_f \lon 0. 
$$
Remark that $R_f \subset {\mathcal E}_N$ is an 
${\mathcal E}_W$-submodule via $f^* : {\mathcal E}_W \to {\mathcal E}_N$.

Now, in $V_f$, 
the ${\mathcal E}_N$-module consisting of vector fields 
along $f$, we have 
$$
VI_f' = \{ v \in V_f \mid i_v\alpha = 0, i_vd\alpha = 0 \}, 
$$
and $VI_f' \subset V_f$ is an ${\mathcal E}_N$-submodule, therefore, 
an ${\mathcal E}_W$-submodule via $f^*$.

\ 

To proceed algebraic calculation, we are going to 
provide also $VI_f$ a module structure. 

As in the previous section, we denote by $X_H$ the contact 
Hamilton vector field with Hamilton function $H$. 

\bep
\label{module} 
$VI_f$ is an ${\mathcal E}_W$-module by the multiplication 
$$
H * v = f^*H\cdot v + (i_v\alpha)(X_H - H\cdot R)\circ f, 
$$
for $H \in {\mathcal E}_W, v \in VI_f$. The multiplication is independent 
on the choice of contact form $\alpha$, but it depends only on 
the contact structure (and on $H, v$). 
Moreover the sequence 
$$
0 \lon VI_f' \lon VI_f \stackrel{e}{\lon} R_f \lon 0
$$
is ${\mathcal E}_W$-exact. 

\enp

\ber
{\rm 
For a constant function 
$c$, we have $X_c = cR$ and $c*v = cv$. 
}
\enr

To verify Proposition \ref{module}, 
we need several lemmas: 

\bel
\label{Lemma}
$i_{H*v}\alpha = f^*H \cdot i_v\alpha$. 
\enl

\Proof 
Since $i_{(HR-X_H)\circ f}\alpha = f^*(i_{HR-X_H}\alpha) = 
f^*(Hi_R\alpha - i_{X_H}\alpha) = f^*(H - H) = 0$, 
we see $i_{H*v}\alpha = i_{f^*H\cdot v}\alpha = 
f^*H\cdot i_v\alpha$. \QED

\bel
\label{Indep}
Set $\alpha' = \lambda\alpha$, for a non-vanishing function $\lambda$. 
Then $i_v\alpha' = f^*\lambda i_v\alpha$ for any vector field along a mapping 
$f : N \to W$. 
If we denote by $R', X'_H$ the Reeb vector field and the contact Hamilton vector 
field of $H$ with respect to $\alpha'$, respectively, and if 
$f : N \to W$ is integral, 
then 
$$
(X'_H - HR')\circ f = \left\{\dfrac{1}{\lambda}(X_H - HR)\right\}\circ f. 
$$
Thefore we have  
$$
(i_v\alpha')(X'_H - HR')\circ f = (i_v\alpha)(X_H - HR)\circ f. 
$$
\enl

\Proof 
That $i_v\alpha' = f^*\lambda i_v\alpha$ follows by Lemma \ref{i and L} (1). 

Set $u = (X'_H - HR')\circ f$ and $v = \left\{\dfrac{1}{\lambda}(X_H - HR)\right\}\circ f$. 
Then, by Lemma \ref{pullback.etc}, 
$i_u\alpha' = i_{X'_H\circ f}\alpha' - i_{(HR')\circ f}\alpha' 
= f^*H - f^*H = 0$. 
Similarly we have $i_v\alpha = 0$. So we have $i_v\alpha' = (f^*\lambda)(i_v\alpha) $. 

We will show $i_ud\alpha' = i_vd\alpha' = f^*(-dH)$. 
Then, since $\alpha'$ is a contact form, we have $u = v$. 

Now in fact, since $f$ is integral, we have $f^*\alpha' = f^*\alpha = 0$,  
and therefore we have, by Lemma \ref{LieReeb}, 
$$
\begin{array}{rcl}
i_ud\alpha' & = & f^*(i_{X'_H - HR'}d\alpha') \\
 & = & f^*(R'(H)\alpha' - dH - Hi_{R'}d\alpha') \\
 & = & f^*(- dH). \\
i_vd\alpha' & = & 
(f^*\lambda)(i_vd\alpha) + i_v(d\lambda\wedge \alpha) \\
 & = & f^*(i_{X_H - HR}d\alpha) + (i_vd\lambda)f^*\alpha 
- (i_v\alpha)f^*(d\lambda) \\
 & = & f^*(R(H)\alpha - dH) \\
 & = & f^*(-dH).
\end{array}
$$

\QED

\ber
{\rm
The terms $(i_v\alpha)X_H\circ f$ and 
$(i_v\alpha) (H\cdot R)\circ f$ do depend on the choice of $\alpha$. 
Just the difference is intrinsically defined as seen in Lemma \ref{Indep}. 
}
\enr

\noindent{\it Proof of Proposition \ref{module}:}
We compare 
$$
(KH)*v = f^*(KH)\cdot v - (i_v\alpha)(KH\cdot R - X_{KH})\circ f
$$
with 
$$
K*(H*v) = f^*K(f^*H\cdot v - (i_v\alpha)(H\cdot R - X_H)\circ f)
- (i_{H*v}\alpha)(K\cdot R - X_H)\circ f. 
$$
By Lemma \ref{Lemma},  
the right hand side of the latter equals to 
$$
f^*(KH)\cdot v - (i_v\alpha)(2KH\cdot R - KX_H - HX_K)\circ f, 
$$
which is equal to 
the right hand side of the former, by Lemma \ref{KH}. 
By Lemma \ref{Lemma}, $e$ is an ${\mathcal E}_W$-epimorphism. 
By Lemma \ref{Indep}, we see the multiplication depends only on the contact 
structure. 
The remaining parts are clear. 
\QED

The following is a consequence of Proposition \ref{module}, 
Lemma \ref{generating function 1} and Proposition \ref{affine function}: 

\bel
If we set 
$$
VH'_{W, f} = \{ X_H \in VH_W \mid H\circ f = 0\}, 
$$
then we have an ${\mathcal E}_W$-exact sequence, 
$$
0 \lon VI_f'/wf(VH'_{W, f}) \lon VI_f/wf(VH_W) \lon R_f/{\mathcal E}_W \lon 0. 
$$
If we set 
$$
VL'_{W, f} = \{ X_H \in VL_W \mid H\circ f = 0\}, 
$$
then we have an ${\mathcal E}_Z$-exact sequence, 
$$
0 \to VI_f'/wf(VH'_{W, f}) \to VI_f/wf(VL_W) 
\to R_f/(E_Z + \sum_{i=1}^n E_Z(p_i\circ f)) \to 0. 
$$
\enl

Let $f : (N, x_0) \to (W, w_0)$ be an integral mapping. 
We define an ${\mathcal E}_W$-homomorphism 
$tf : V_N \to VI_f$ by $tf(\xi) := f_*(\xi), \xi \in V_N$. 

\bel
Let $f : (N, x_0) \to (W, w_0)$ be an integral map-germ. 
Then $tf(V_N) \subseteq VI'_f$. 
\enl

\Proof
Take $f_*(\xi) \in tf(V_N)$. Then we have 
$e(f_*(\xi)) = i_{f_*(\xi)}\alpha = i_{\xi}f^*\alpha = 0$. 
\QED

Under a condition, the converse inclusion holds: 

\bep 
\label{codim 2}
Let $f : (N, x_0) \to (W, w_0)$ be an integral map-germ. 
Suppose that 
$f$ is diffeomorphic to an analytic map-germ $f' : (\KK^n, 0) 
\to (\KK^{2n+1}, 0)$ {\rm (}not necessarily integral{\rm )} 
such that the codimension of the 
singular locus of the complexification $f'_{\C}$ of $f'$ is greater than 
or equal to $2$. 
Then we have $VI'_f \subseteq tf(V_N)$. 
Therefore we have an isomorphism of ${\mathcal E}_W$-modules
$$
VI_f/\{ tf(V_N) + wf(VH_W)\} \cong R_f/{\mathcal E}_W, 
$$
and an isomorphism of ${\mathcal E}_Z$-modules
$$
VI_f/\{ tf(V_N) + wf(VL_W)\} \cong R_f/(E_Z + \sum_{i=1}^n E_Z(p_i\circ f)). 
$$
\enp

\Proof
Let $v \in VI'_{f}$. 
Set 
$$
v = (p\circ f, q\circ f, r\circ f; \phi\circ v, \xi\circ v, s\circ v). 
$$
Then, since $e(v) = s\circ v - \sum(p\circ f)(\xi\circ v) = 0$, we have 
$$
\sum_{i=1}^n(\xi_i\circ v)d(p_i\circ f) - \sum_{i=1}^n(\phi_i\circ v)d(q_i\circ f) = 0. 
$$ 
This means, for any regular point $x \in (\KK^n, 0)$, that $v(x) \in D_{f(x)}$ 
and $v(x)$ belongs to the skew orthogonal complement to 
$f_*(T_x\KK^n)$ 
with respect to the symplectic structure $\sum_{i=1}^n dp_i \wedge dq_i$ on $D$. 
Therefore we have $v(x) \in f_*(T_x\KK^n)$. 
Since $f$ and $f'$ are diffeomorphic, any vector field 
in $VI'_f$ is transformed to a vector field along $f'$ which is tangent to 
the image of $f'$ off the singular locus of $f'$. 

Let $v \in V_{f'_{\C}}$. This means 
that $v : (\C^n, 0) 
\to T\C^{2n+1}$ is a holomorphic vector field along  
$f'_{\C} : (\C^n, 0) \to (\C^{2n+1}, 0)$. 
Suppose, for each regular point $x \in (\C^n, 0)$ of $f'_{\C}$ 
that $v(x) \in f'_{\C *}(T_x\C^n)$. 
Then we can find a vector field $w$ over 
$\C^n \setminus {\mbox {\rm Sing}}(f'_{\C})$ 
satisfying $v = (f'_{\C})_*(w)$ on $\C^n \setminus {\mbox {\rm Sing}}(f'_{\C})$, 
where ${\mbox {\rm Sing}}(f'_{\C})$ is the locus of singular points of $f'_{\C}$. 
Since ${\mbox {\rm Sing}}(f'_{\C})$ is of codimension $\geq 2$ in $\C^n$, 
$w$ extends to a holomorphic vector field on $(\C^n, 0)$ still called $w$, 
by Hartogs theorem. Then we have $v = f'_{\C *}(w)$. 
This proves that $VI'_f \subseteq tf(V_N)$ in the case $\KK = \C$. 

In the case $\KK = \R$, we set  
$T \subset V_{f'}$ as the set of vector fields along $f'$ such that, 
for each regular point $x \in (\R^n, 0)$ 
of $f'$, $v(x) \in f'_*(T_x\R^n)$. 

Take $v \in T$. Suppose $v$ is real analytic. 
Then considering 
the complexification of $v$, we see that there exists a real analytic 
$w \in V_n$ such that $v = f'_*(w)$ over $(\R^n, 0)$. 
This means that $T$ is generated formally by $tf'(V_n)$ in the sense 
of \cite{Malgrange}, and, by Whitney's spectral theorem,  we have that 
$T$ is contained in the closure of $tf'(V_n)$ for a representative of $f'$. 
Since $tf'(V_n)$ itself is closed, we see $T \subseteq tf'(V_n)$. 
This shows that $VI_f \subseteq tf(V_N)$. 

The remaining parts are clear. 
\QED

\ 

We call $f$ {\it infinitesimally contact stable} if 
$$
VI_f = tf(V_N) + wf(VH_W). 
$$
Then we have: 

\bec
\label{prop codim 2}
Let $f : (N, x_0) \to (W, w_0)$ be an integral mapping. 
Then the condition {\rm (}ca{\rm )} of Proposition \ref{contact stable}
implies that $f$ is infinitesimally contact stable, namely the condition {\rm (}ics{\rm )}. 
\enc

\Proof
Since $R_f = f^*{\mathrm E}_W$, we see  
$0 = R_f/{\mathcal E}_W \cong VI_f/\{ tf(V_N) + wf(VH_W)\}$. 
Therefore we have $VI_f = tf(V_N) + wf(VH_W)$. 
\QED

\bel
If an integral map-germ of corank at most one
$f : (N, x_0) \to W$ is infinitesimally contact stable 
then $f$ is a finite map-germ. 
\label{Finitefin}
\enl

\Proof
By taking generating functions of both sides of the equality 
$VI_f = tf(V_N) + wf(VH_W)$,  
we have $R_f = f^*{\mathcal E}_W$. Therefore $R_f$ is a 
finite ${\mathcal E}_W$-module. Therefore, by Lemma 
\ref{R_f finite}, we see $f$ is finite. 
\QED

\ 

Let 
$(f_t)$ be an integral deformation of $f$. 
To show $f$ is homotopically contact (resp. Legendre) stable, 
we need to find a deformation $(\sigma_t)$ of $\id_N$ 
and an integral deformation $(\tau_t)$ of $\id_W$ 
(resp. an integral deformation $\tau_t$ of $\id_W$ 
covering a deformation $(\bar{\tau}_t)$ of $\id_Z$ via 
$\pi : W \to Z$) satisfying 
$\tau_t^{-1}\circ f_t\circ\sigma_t = f$. 
For this, it is sufficient to solve 
$df_t/dt = 
\eta_t\circ f_t - Tf_t\circ\xi_t 
( = wf_t(\eta_t) - tf_t(\xi_t) ) 
: N\times\KK \to TW$ 
with 
$
\xi_t \in V_N$ 
and 
$\eta_t \in VH_W$  
(resp. 
$\eta_t \in VL_{W}$), 
(cf. \cite{Mather}). 

\ 

For an unfolding $F = (f_t, t)
 : N\times J \to W\times J$, \, $t \in J = (\KK,0)$, we set 
\[
VI_{F/J} = \{v : N\times J \to TW \  \mid \ v_t \in VI_{f_t}, \, t \in J\}.
\] 
If 
$(f_t)$ is an integral deformation of $f$, 
then 
we have 
$(df_t/dt)_{t\in J} \in VI_{F/J}$. 
We define an ${\mathcal E}_{W\times J}$-module structure on 
$VI_{F/J}$ by 
$$
a_t*v_t = (f_t^*a_t)\cdot v_t + (i_{v_t}\alpha)(X_{a_t} - a_t\cdot R)\circ f_t, 
$$ 
for $v_t \in VI_{F/J}, a_t \in {\mathcal E}_{W\times J}$. 
Compare with Proposition \ref{module}. 
Then we have 
\bec
If $f$ is finite and of corank at most one, then the 
quotient 
$VI_{F/J}$ is a finite ${\mathcal E}_{W\times J}$-module. 
\label{Fin}
\enc

Now assume $f$ is integral and 
$f_t$ is an integral deformation of $f$. 
We define $tF/J : V_N \to VI_{F/J}$ by $v \mapsto 
(tf_t(v))_{t \in J}$. 
We set 
$$
S_{F/J} = VI_{F/J}/((wF/J)(VH_W) + (tF/J)(V_N)), 
$$ 
which is an 
${\mathcal E}_{W\times J}$-module, and set 
$$
S_f = VI_f/(wf(VH_W) + tf(V_N)),
$$ 
which is an ${\mathcal E}_W$-module. 
Then we have: 
\bel
The quotient 
$S_{F/J}/m_JS_{F/J}$ is isomorphic to $S_f$ 
as an ${\mathcal E}_W$-modules.
\label{Quotient}
\enl 
\Proof
Consider the morphism 
$\Phi : 
S_{F/J}
\to S_f
$
defined by 
$\Phi([v_t]) = [v_t\vert_{t=0}]$. 
We will show that the kernel of $\Phi$ is equal to 
$m_JS_{F/J}$. 
Let $v_t \in 
VI_{F/J}
$. 
Assume $v_t\vert_{t=0} = wf(\eta) + tf(\xi)$, 
for some 
$\xi \in V_N$, $\eta \in VI_W$. 
Set 
$w_t = v_t - wf_t(\eta) - tf_t(\xi)$. Then 
$w_t\vert_{t=0} = 0$. 
Therefore $w_t = tw_t'$, for some 
$w_t' \in V_{f_t}$. 
We see 
$\Pi_{\sharp}(w_t') \in VI_{g_t}$. 
Here $g_t = \Pi\circ f_t$ is the family of isotropic map-germs induced 
from $f_t$. 
In fact 
$\Pi_{\sharp}(w_t)^\flat = {t\Pi_{\sharp}(w_t')}^\flat$ and so 
$0 = 
(\Pi_{\sharp}(w_t)^\flat)^*d\theta_{T^*Q} 
= t({\Pi_{\sharp}(w_t')}^\flat)^*d\theta_{T^*Q}$. Thus 
$({\Pi_{\sharp}(w_t')}^\flat)^*d\theta_M = 0$. This means 
$w_t' \in VI_{f_t}$. 
Since $x$-derivative of $t$ is equal to zero, 
we have 
$w_t = tw_t = t*w_t'$ and 
$[v_t] = [w_t] = t[w_t'] \in m_JS_{F/J}$. 
\QED

\section{Relation to Isotropic Mappings.} 
\label{Isotropic map}

Let $Q$ be a manifold of dimension $n$. Then $T^*Q\times \KK \cong J^1(Q, \KK) 
\subset PT^*(Q\times \KK)$ has the canonical 
contact structure, whereas $T^*Q$ has the canonical symplectic 
structure $\omega = d\theta_Q$, $\theta_Q$ being Liouville form on $Q$, 
$\theta_Q = \sum_{i=1}^n p_idq_i$, for a system of  local symplectic coordinates. 
A contact form on $T^*Q\times \KK$ is given by $dr - \theta_Q$, 
for the coordinate $r$ on $\KK$.  

Let $g : N \to T^*Q$ be a mapping from a manifold $N$ of 
dimension $n$. Then $g$ is called {\it isotropic} 
if $g^*\omega = 0$. The singularities of isotropic 
mappings of corank at most one is studied in \cite{Inv} in detail. 
In particular, we have a series of singularities, 
\lq\lq open Whitney umbrellas", which are symplectic counterparts 
of objects we have introduced in this paper. 

Two isotropic map-germs $g : (N, x_0) 
\to T^*Q$ and $g' : (N, x'_0) \to T^*Q$ 
are called {\it symplectomorphic} (or {\it symplectically  
equivalent}) 
if there exist a symplectomorphism 
$$\tau : (T^*Q, g(x_0)) 
\to (T^*Q, g'(x'_0))$$ and a diffeomorphism 
$\sigma : (N, x_0) \to (N, x'_0)$ satisfying 
$\tau\circ f = f'\circ \sigma$. 
Then we call also the pair $(\sigma, \tau)$ a 
symplectomorphism between $g$ and $g'$.

Let $f : (N, x_0) \to T^*Q\times \KK$ be a map-germ. 
Set $g : (N, x_0) \to T^*Q$ to be $g = \Pi\circ f$, where 
$\Pi : T^*Q\times \KK \to T^*Q$ is the natural projection along 
the flow of Reeb vector field $\dfrac{\pa}{\pa r}$.

Then we have by \cite{Inv}: 

\bel
\label{int-isot}
\par 
{\mbox{\rm (1)}} 
$f$ is an integral map-germ if and only if 
$g$ is an isotropic map-germ.

{\mbox{\rm (2)}}
If $g = \Pi\circ f$ and $g' = \Pi\circ f'$ 
are symplectomorphic, then $f$ and $f'$ are contactomorphic. 

{\mbox{\rm (3)}}
$R_f = R_g$.

{\mbox{\rm (4)}}
$f$ is an open Whitney umbrella of type $k$ {\rm (}as an integral map-germ{\rm )} 
if and only if $g$ is an open Whitney umbrella of type $k$ 
{\rm (}as an isotropic map-germ{\rm )}. 
In particular, 
$f$ is a Legendre immersion if and only if $g$ is a Lagrange immersion.

\enl

\ber
{\rm 
The converse of (2) of Lemma \ref{int-isot} does not hold in general. 
For example, consider integral map-germs  
$f_\lambda : (\KK, 0) \to (\KK^3, 0), \lambda > 0$ 
defined by $g(t) = (t^3, t^7 + \lambda t^8, \dfrac{3}{10}t^{10} + 
\dfrac{3}{11}\lambda t^{11})$. Then $g_\lambda = \Pi\circ f_\lambda : 
(\KK, 0) \to (\KK^2, 0), $ 
$g_\lambda(t)  = (t^3, t^7 + \lambda t^8)$  
is not symplectomorphic to $g_{\lambda'}$ if 
$\lambda' \not= \lambda$, while all $f_\lambda$ 
are contactomorphic to each other (\cite{IJ}\cite{CSLCC}). 
}
\enr

Set $W = T^*Q\times \KK$. 
The projection $\Pi : W \to T^*Q$ induces the projection 
$\Pi_* : TW \to T(T^*Q)$; by using local coordinates, 
it is given by 
$$
\Pi_*(p, q, r; \phi, \xi, s) = (p, q; \phi, \xi). 
$$
Then $\Pi_*$ induces $\KK$-linear mapping 
$\Pi_{\sharp} : V_f \to V_g$ by $\Pi_{\sharp}(v) =  \Pi_*\circ v, (v \in V_f)$. 

Now we observe the following: 

\bel
\label{Diagram} 
$\Pi_{\sharp}$ restricts to a $\KK$-linear epimorphism 
$\Pi_{\sharp} : VI_f \to VI_g$, to an ${\mathcal E}_N$- isomorphism 
$\Pi_{\sharp} : VI'_f \to VI'_g$ and 
${\mathcal E}_{T^*Q}$-epimorphism 
$\Pi_{\sharp} : VI_f \to VI_g/wg(VH_{T^*Q})$ 
over the ring morphism 
$\Pi^* : {\mathcal E}_{T^*Q} \to {\mathcal E}_W$. 
Furthermore we have the following commutative 
diagram which consists of exact sequences: 
{\small 
$$
\begin{array}{ccccccccc}
 & & 0 & & 0 & & 0 & & \\
 & & \downarrow & & \downarrow & & \downarrow & & \\
0 & \rightarrow & wf(VH_{W, g}') & \rightarrow 
& wf(VH_W') & \to & f^*{\mathcal E}_W 
& \rightarrow & 0 \\
 & & \downarrow & & \downarrow & & \downarrow & & \\ 
0 & \rightarrow & VI'_f & \rightarrow & VI_f & \stackrel{\widetilde{e}}{\to} & R_f & \rightarrow & 0 \\
  &  & \mapdown{ } &  & \mapdown{\overline{\Pi}_{\sharp}} &  & \mapdown{ } & &  \\
0 & \rightarrow & VI'_g/wg(VH'_{T^*Q, g}) & \rightarrow
& VI_g/wg(VH_{T^*Q})
& \stackrel{\overline{e}}{\to} & R_g/g^*{\mathcal E}_{T^*Q} & \rightarrow & 0 \\
 & & \downarrow & & \downarrow & & \downarrow & & \\
 & & 0 &  & 0 & & 0 & & \\
\end{array}
$$
}

The kernel of $\Pi_{\sharp}$ is generated by $R\circ f = \dfrac{\pa}{\pa r}\circ f$ 
over $\R$. 
\enl

\Proof
We show that $\Ker(\Pi_{\sharp}) = \left\langle \dfrac{\pa}{\pa r}\circ f\right\rangle_{\R}$. 
Let $v  = (p\circ f, q\circ f, r\circ f; \phi, \xi, s) \in VI_f$. 
Recall that $d(s - (p\circ f)\xi) + \xi d(p\circ f) - \phi d(q\circ f) = 0$. 
Suppose that $\Pi_*\circ v = 0$. Then $\xi = 0, \phi = 0$. 
Then we have $ds = 0$. Thus $s$ is constant. The remaining parts are clear. 
\QED

We have also 

\bel
\label{Lifting}
For any $\eta \in VH_{T^*Q}$ {\rm (}resp. $\eta \in 
VL_{T^*Q}${\rm )}, there exists an $\widetilde{\eta} \in VH_W$ 
{\rm (}resp. $\widetilde{\eta} \in VL_W${\rm )}, 
such that 
$\Pi_{\sharp}wf(\widetilde{\eta}) = wg(\eta)$. Here 
$wf(\widetilde{\eta}) = \widetilde{\eta}\circ f$ and 
$wg(\eta) = \eta\circ g$. $VL_{T^*Q}$ means the set of 
Lagrange vector fields of the Lagrange fibration $T^*Q \to Q$ {\rm (\cite{Inv})}. 

If $\eta \in m_{T^*Q}*VH_{T^*Q}$ {\rm (}resp. $\eta \in 
m_Q*VL_{T^*Q}${\rm )}, then we can take $\widetilde{\eta}$ from 
$m_W^2* VH_W$ {\rm (}resp. from $m_Z^2* VL_W${\rm )}, where $Z = Q\times \R$. 
\enl

\Proof 
If $\eta$ has a symplectic Hamiltonian function $H$ on $T^*Q$, 
$H(0) = 0$, 
then we may set $\widetilde{\eta} = X_{\Pi^*H}$, 
the contact Hamiltonian vector field for the pull-back $\Pi^*H$ 
of $H$ by $\Pi$. 
\QED

\bel
\label{symp-cont}
Let $f : (N, x_0) \to W = T^*Q \times \KK$ be an integral mapping. 
If $g = \Pi\circ f : (N, x_0) \to T^*Q$ is infinitesimally symplectically 
{\rm (}resp. Lagrange{\rm )} 
stable, then $f$ is infinitesimally contact {\rm (}resp. Legendre{\rm )} stable. 
\enl

\Proof
Suppose $g$ is infinitesimally symplectic (resp. Lagrange) stable. 
Then,  for any $v \in VI_f$, there exist $\xi \in V_N$ and 
$\eta \in VH_{T^*Q}$ (resp. $\eta \in VL_{T^*Q}$) 
satisfying $\Pi_{\sharp}v = tg(\xi) + wg(\eta)$. 
Then we see, using Lemma \ref{Lifting}
$$
\Pi_{\sharp}(v - tf(\xi) - wf(\widetilde{\eta})) = 0. 
$$
Thefore, by Lemma \ref{Diagram}, there exists $s_0 \in \R$ such that 
$$
v - tf(\xi) - wf(\widetilde{\eta}) = s_0\dfrac{\pa}{\pa r}\circ f. 
$$
Thus we have $v = tf(\xi) + wf(\widetilde{\eta} + s_0\dfrac{\pa}{\pa r})$. 
\QED

\bep
\label{prop OWU}
If $f$ is an open Whitney umbrella, then 
we have $R_f = f^*{\mathcal E}_W$. Therefore 
$f$ satisfies 
the condition {\rm (}ca{\rm )} of Proposition \ref{contact stable}. 
\enp

\Proof
$R_f = R_g = g^*{\mathcal E}_Z \subseteq f^*{\mathcal E}_W \subseteq R_f$. 
\QED

\bec
If $f$ is an open Whitney umbrella, then 
$f$ is infinitesimally contact stable.  
Moreover we have the isomorphism 
$$
VI_f/(tf(V_N) + wf(VL_W)) \cong R_f/(E_Z + \sum_{i=1}^n E_Z(p_i\circ f))
$$
of ${\mathcal E}_Z$-modules 
via $\pi^* : {\mathcal E}_Z \to {\mathcal E}_W$. 
\enc

\section{Integral Jets.}
\label{Integral jet} 

We consider the {\it integral jet space}: 
$$
J^r_I(n, 2n+1) := \{j^rf(0) \mid 
f : (\KK^n, 0) \to (\KK^{2n+1}, 0) {\mbox{\rm \ integral \ of \ corank}}\leq 1\}. 
$$ 
Then $J^r_I(n, 2n+1)$ is a submanifold of $J^r(n, 2n+1)$. 

\ber
{\rm 
The projection $\Pi^r : 
J^r(n, 2n+1) \to J^r(n, 2n)$ 
defined by $\Phi^r(j^rf(0)) := j^r(\Pi\circ f)(0)$ 
induces a diffeomorphism of 
$J^r_I(n, 2n+1)$ and the isotropic jet space $J^r_I(n, 2n) \subset J^r(n, 2n)$ 
(\cite{Banach}). In fact, for any $j^rg(0) \in J^r_I(n, 2n)$, 
we set $j^rf(0) = j^r(g, e)(0)$, where $e$ is the generating function 
of $g$, $de = g^*\theta_Q, e(0) = 0$. 
Then $j^rf(0) \in J^r_I(n, 2n+1)$ and 
$\Pi^r(j^rf(0)) = j^rg(0)$. 
}
\enr

Let 
$f : (N, x_0) \to (W, w_0)$ be an integral map-germ of corank at most one. 
Then we set 
$$
VI_f^s = \{ v \in VI_f \mid j^sv(x_0) = 0\} = VI_f \cap m_N^{s+1}V_f, \ 
(s = 0, 1, 2, \dots). 
$$

Let $z = j^rf(x_0) \in J^r_I(n, 2n+1)$. 
Define $\pi_r : VI_f^0 \to T_zJ^r(n, 2n+1)$ as follows: 
For each $v \in VI_f^0$, take an integral deformation $(f_t)$ of $f$ with 
$v = \left.\dfrac{df_t}{dt}\right\vert_{t=0}$, 
and set $\pi_r(v) = \left.\dfrac{d(j^rf_t(x_0))}{dt}\right\vert_{t=0}$. 
Then the image of the linear map $\pi_r$ coincides with 
$T_zJ^r_I(n, 2n+1)$. 

Let $z \in J_I^r(n, 2n+1)$ and $z = j^rf(x_0)$ for a $f : 
(N, x_0) \to (W, w_0)$. Then 
under the identification $T_zJ^r(n, 2n+1) \cong m_NV_f/m_N^{r+1}V_f$ we have 
$$
T_zJ_I^r(n, 2n+1) \cong VI^0_f/VI^r_f. 
$$
If we denote by ${\mathcal C}^rz$ (resp, ${\mathcal{L}}^rz$) the orbit of $z$ 
under the contactomorphisms (resp. Legendre diffeomorphisms), we have 
$$
T_z{\mathcal C}^rz \cong \{(tf(m_NV_N) + wf(m_W^{(2)}*VH_W)) + VI^r_f\}/VI^r_f, 
$$
$$
T_z{\mathcal{L}}^rz \cong \{(tf(m_NV_N) + wf(m_Z^{(2)}*VL_W)) + VI^r_f\}/VI^r_f. 
$$

Set $z = j^rf(x_0)$. 
For $(w, v) \in T_{x_0}N \oplus VI_f$, take a curve $x_t$ in $N$ with the 
velocity vector $w$ at $t = 0$ and take an integral deformation 
$f_t$ of $f$ with $v = \left.\dfrac{df_t}{dt}\right\vert_{t=0}$ (cf. 
\cite{Inv}, Lemma 3.4), 
and define a linear map 
$$
\Pi_r : T_{0}N \oplus VI_f \to T_zJ^r(N, W), 
$$
by 
$$
\Pi_r(w, v) = \dfrac{j^rdf_t(x_t)}{dt}\vert_{t=0}. 
$$
Then $\Pi_r(T_{0}N \oplus VI_f) = T_zJ^r_I(N, W)$ and 
$\Ker \Pi_r = \{ 0\}\oplus VI_f^r$. 
Moreover we have, for the Legendre equivalence class,   
$$
[z] = 
\{ j^rf'(x) \mid 
x \in N, \ f' {\mbox{\rm \ is Legendre equivalent to\  }} f \}
$$ 
in $J_I^r(N, W)$, 
$$
T_z[z] = \Pi_r(T_{x_0}N \oplus(tf(m_NV_N) + wf(VL_W))). 
$$
For the jet extension 
$j^rf : (N, x_0) \to J^r_I(N, W)$, we have 
$$
(j^rf)_*(\dfrac{\pa}{\pa x_i}) = 
\Pi_r(\dfrac{\pa}{\pa x_i}, f_*(\dfrac{\pa}{\pa x_i})). 
$$

\bel
\label{transversality cond}
The transversality condition {\rm(t}$_r${\rm )} is equivalent to 
the condition 
$$
VI_f = tf(V_N) + wf(VL_W) + VI_f^r. 
$$ 
\enl

\Proof 
The condition (t$_r$) that $j^rf$ is transverse to $[z] = [j^rf(x_0)]$ 
at $x_0$ is equivalent to the condition 
$$
(j^rf)_*(T_{x_0}N) + T_z[z] = T_zJ^r_I(N, W), 
$$
and to the condition  that 
$$
(\Pi_r)^{-1}((j^rf)_*(T_{x_0}N)) + 
T_{x_0}N\oplus(tf(m_NV_N) + wf(VL_W)) 
+ \{ 0\}\oplus VI_f^r
$$
coincides with $T_{x_0}N\oplus VI_f$. This condition is equivalent to 
that 
$$
VI_f = \langle f_*(\dfrac{\pa}{\pa x_1}), \dots, 
f_*(\dfrac{\pa}{\pa x_n})\rangle_{\KK} 
+ tf(m_NV_N) + wf(VL_W) + VI_f^r, 
$$
namely that 
$$
VI_f = tf(V_N) + wf(VL_W) + VI_f^r. 
$$
\QED

Similarly we have the following: 
\bel
\label{cont-trans} 
The condition that $j^rf$ is transverse to the contactomorphism-orbit through 
$j^rf(x_0)$ at $x_0$ 
is equivalent to the condition 
$$
VI_f = tf(V_N) + wf(VH_W) + VI_f^r. 
$$
\enl

Moreover the transversality condition on $j^rf$ 
implies that $f$ is an open Whitney umbrella: 

\bep
\label{owu,ct}
Let $f : (N, x_0) \to (W, w_0)$ 
be an integral map-germ of corank $\leq 1$, 
and $k$ a non-negative integer. 
If the $k+1$ extension 
$j^{k+1}f : (N, x_0) \to 
J^{k+1}_I(N, W)$ is transverse 
to the contactomorphism-orbit through 
$j^{k+1}f(x_0)$, then 
$f$ is an open Whitney umbrella of type $\leq k$. 
\enp

\Proof 
Since $f$ is an integral map-germ of corank $\leq 1$, 
$f$ is contactomorphic to 
$f' : (\KK^n, 0) \to (\KK^{2n+1}, 0)$ 
with 
$$
\varphi := (q_1, \dots, q_{n-1}, q_n, p_n)\circ f' = (x_1, \dots, x_{n-1}, u(x), v(x)). 
$$
Since $f'$ is contactomorphic to $f$, 
also $j^{k+1}f' : (\KK^n, 0) \to 
J^{k+1}_I(\KK^n, \KK^{2n+1})$ is transverse 
to the contactomorphism-orbit through 
$j^{k+1}f'(0)$, therefore 
to ${\mathcal K}$-orbit through 
$j^{k+1}f'(0)$. 
Then we see 
$j^{k+1}\varphi : (\KK^n, 0) \to 
J^{k+1}(\KK^n, \KK^{n+1})$ is transverse 
to ${\mathcal K}$-orbit through 
$j^{k+1}\varphi(0)$. 
Then $f'$ is an open Whitney umbrella of type $\leq k$. 
Therefore $f$ is an open Whitney umbrella of type $\leq k$. 
\QED

For an $n$-dimensional manifold $N$ and a contact manifold $W$ of dimension $2n+1$, 
we set
$$
C^\infty_I(N, W)^1 := 
\{ f : N \to W \mid f {\rm \ is \ integral \ of \ corank } \leq 1\}. 
$$
We endow $C^\infty_I(N, W)^1$ with the relative topology of the 
Whitney $C^\infty$ topology of $C^\infty(N, W)$. 
Then we have the following Legendre transversality theorem: 

\bep
Let $r$ be a non-negative integer and $U$ a locally finite family of submanifolds 
of $J^r(N, W)$. Then 
$$
T_U := \{ f \in C^\infty_I(N, W)^1 \mid j^rf {\rm 
\ is \ transverse \ to \ all \ of \ } U \} 
$$
is dense in $C^\infty_I(N, W)^1$. 
\enp

\section{Finite determinacy.}
\label{Finite determinacy}

\bel
\label{k+1} 
Let $f, f' : (N, x_0) 
\to W$ be integral map-germs. 
If $f$ is an open Whitney umbrella of type $k$ and 
$j^{k+1}f'(x_0) = j^{k+1}f(x_0)$, 
then $f'$ is an open Whitney umbrella of type $k$. 
\enl

\Proof
By definition there exist a contactomorphism $(\sigma, \tau)$ 
such that $\tau\circ f\circ\sigma^{-1} = 
f_{n,k}$, the normal form. 
Set $f'' = \tau\circ f'\circ\sigma^{-1}$. Then 
$j^{k+1}f''(x_0) = j^{k+1}f_{n,k}(x_0)$. 
Set $\varphi = (q_1, \dots, q_{n-1}, q_n, p_n)\circ f'' : (N, x_0) \to 
\KK^{n+1}$ and 
$\varphi_{n,k} = (q_1, \dots, q_{n-1}, q_n, p_n)\circ f_{n,k}$. 
Then $j^{k+1}\varphi(x_0) = j^{k+1}\varphi_{n,k}(x_0)$. 
Now $j^k\varphi : (N, x_0) 
\to J^k(N, \KK^{n+1})$ 
is transverse at $x_0$ to 
Thom-Boardman strata as well as $j^k\varphi_{n, k}$ is. 
In \cite{Camb}, we have shown that $g'' = (q, p)\circ f'' : 
(N, x_0) \to \KK^{2n}$ 
is symplectomorphic to $g_{n, k} = (q, p)\circ f_{n,k}$. 
Then $f''$ and $f_{n,k}$ are contactomorphic. Since 
$f$ and $f''$ are contactomorphic, we see $f$ and $f_{n,k}$ are contactomorphic, 
therefore $f$ is an open Whitney umbrella of type $k$. 
\QED

\ 


An integral map-germ $f : (N, x_0) \to (W, w_0)$ 
is called {\it $r$-determined by contactomorphisms} 
if, for any integral map-germ $f' : (N, x_0) \to (W, w_0)$ 
with $j^rf'(x_0) = j^rf(x_0)$, $f$ and $f'$ 
are contactomorphic. 

Let $\pi : (W, w_0) \to (Z, z_0)$ be a fixed 
Legendre fibration. 
An integral map-germ $f : (N, x_0) \to (W, w_0)$ 
is called {\it Legendre $r$-determined} if, 
for any integral map-germ $f' : (N, x_0) \to (W, y_0)$ 
with $j^rf'(x_0) = j^rf(x_0)$, then $(f', \pi)$ and $(f, \pi)$ 
are Legendre equivalent. 

Then we have: 
\bel
An open Whitney umbrella of type $k$ is $(k+1)$-determined by contactomorphisms. 
\enl

\Proof
Suppose $f$ is an open Whitney umbrella of type $k$. 
Let $f' : (N, x_0) \to W$ be an integral map-germ with 
$j^{k+1}f'(x_0) = j^{k+1}f(x_0)$. Then $f'$ is also an open Whitney umbrella of type $k$. 
Therefore both $f$ and $f'$ are contactomorphic to the normal form $f_{n,k}$. 
Thus $f$ and $f'$ are contactomorphic. 
\QED

\bel
\label{r-determined}
Let $f : (N, x_0) \to (W, y_0)$ be an open Whitney umbrella. 
Suppose that $f$ is infinitesimally Legendre stable, namely that 
$$
VI_f = tf(V_N) + wf(VL_W). 
$$
Take a positive integer $r$ satisfying 
$$
f^*{\mathcal E}_W \cap m_N^{r+1} \subseteq 
f^*m_W^{n+2}. 
$$
Then we have

{\rm (1)} $R_f = f^*{\mathcal E}_W$ is generated as ${\mathcal E}_Z$-module 
by $1, p_1\circ f, \dots, p_n\circ f$. 

{\rm (2)} $m_W^{n+1}R_f \subseteq m_ZR_f$. 

{\rm (3)}
$ 
VI_f^r \subseteq tf(m_NV_N) + wf(VL_W \cap m_W^{(2)}*VH_W). 
$

{\rm (4)} $f$ is Legendre $r$-determined. 

\enl

\Proof 
(1) : 
Taking generating functions both sides of $VI_f = tf(V_N) + wf(VL_W)$, 
we have 
$$
R_f = \langle 1, p_1\circ f, \dots, p_n\circ f \rangle_{{\mathcal E}_Z}. 
$$
Moreover, since $f$ is an open Whitney umbrella, we have $R_f = f^*{\mathcal E}_W$ 
(Lemma \ref{prop OWU}). 

(2) : 
Set $Q_f := R_f/m_ZR_f$. Then $Q_f$ is generated by 
$1, p_1\circ f, \dots, p_n\circ f$ over $\KK$. Therefore 
$\dim_{\KK}Q_f \leq n+1$. 
Considering the sequence
$$
Q_f \supseteq m_WQ_f \supseteq \cdots \supseteq m_W^{n+1}Q_f, 
$$
and using Nakayama's lemma, 
we have $m_W^{n+1}Q_f = 0$. Therefore we have 
$m_W^{n+1}R_f \subseteq m_ZR_f$. 

(3) : Let $v \in VI_f^r$. Then the generating function $e(v) = i_v\alpha$ 
of $v$ belongs to $f^*{\mathcal E}_W \cap m_N^{r+1}$. 
Now 
$$
f^*{\mathcal E}_W \cap m_N^{r+1} \subseteq 
f^*(m_W^{n+2}) \subseteq m_Zf^*{m_W}. 
$$
Therefore there exist functions $a_1, \dots, a_s \in m_Z$ 
and $b_1, \dots, b_s \in m_W$ such that 
$$
e(v) = (a_1b_1 + \dots + a_sb_s)\circ f. 
$$
For each $b_j\circ f$, there exist $c_{j0},  
c_{j1}, \dots, c_{jn} \in 
{\mathcal E}_Z$ satisfying 
$$
b_j\circ f = c_{j0}\cdot 1 + c_{j1}(p_1\circ f) + \cdots + c_{jn}(p_n\circ f). 
$$
Note that, since $b_j(x_0) = 0$, we see $c_{j0}(x_0) = 0$, therefore 
$c_{j0} m_Z$. 
Set 
$$
h = \sum_{j=1}^s a_j(c_{j0} + c_{j1}p_1 + \cdots + c_{jn}p_n). 
$$
Then $h$ is an affine function with respect to $p_1, \dots, p_n$ 
and $h \in m_W^2$. 
So the Hamilton vector field $X_h$ belongs to 
$VL_W \cap m_W^2*VH_W \subseteq VL_W \cap m_W^{(2)}*VH_W$. 
Then the generating function of $u := v - X_h\circ f$ is equal to zero. 
Then the vector field $u$ is tangent to $f$ along the regular locus 
of $f$. Since $f$ is an open Whitney umbrella, $f$ is analytic and 
the singular locus of the complexification of $f$ is at least $2$. 
Therefore there exists a vector field $\xi \in V_N$ 
satisfying $u = tf(\xi)$ (Proposition \ref{codim 2}). 
So we have $v = tf(\xi) + wf(X_h)$. 
Remark that, since $f$ is an open Whitney umbrella, the kernel of 
the differential mapping $f_* : T_{x_0}N \to T_{w_0}W$ 
is not tangent to 
the Boardman strata containing $x_0$. 
The vector $\xi(x_0)$ belongs to the kernel. On the other hand 
$\xi(x_0)$ must be tangent to the Boardman stratum. 
Therefore we have $\xi(x_0) = 0$ namely 
$\xi \in m_NV_N$. 
Thus we have 
$$
v \in tf(m_NV_N) + wf(VL_W \cap m_W^{(2)}*VH_W). 
$$

(4) : 
Let $f' : (N, x_0) \to W$ be an integral map-germ with 
$j^rf'(x_0) = j^rf(x_0)$. Note that $r \geq n + 1$. 
Therefore $f'$ is also an open Whitney umbrella of the same type as $f$. 
By the argument of Proposition 3.5 in \cite{Banach2}, 
we can connect $f$ and $f'$ by a family of integral map-germs 
$f_t$ satisfying 
$$
VI_{f_t}^r \subseteq tf_t(m_NV_N) + wf_t(VL_Z \cap m_W^{(2)}*VH_W). 
$$
Thus by the homotopy method we see $f$ and $f'$ are Legendre equivalent. 

\section{Contact and Legendre stabilities.} 
\label{Legendre}

First we show the following: 

\bel
\label{stableOWU}
Let $f : (N^n, x_0) \to W^{2n+1}$ be an integral map-germ of 
corank at most one. If 
$f$ is contact stable then $f$ is an open Whitney umbrella.  
\enl

\Proof
Because all notions involved are local and invariant under the contactomorphisms, 
we may assume, by the Darboux theorem, 
 $f : (\KK^n, 0) \to (\KK^{2n+1}, 0)$, $f^*\alpha = 0$, 
$\alpha = dr - pdq$ and $f$ is of corank $\leq 1$. 
Take a representative $f : U \to \KK^{2n+1}$ of the germ $f$. 
We may assume the representative is also integral and of corank $\leq 1$.  

Set $g = (p\circ f, q\circ f) : U \to T^*\KK^n$. Then 
$g$ is isotropic and of corank $\leq 1$. 
Here $g$ is called {\it isotropic} if $g^*\omega = 0$, for the 
symplectic form $\omega = d(pdq)$ on $T^*\KK^n$. 
In fact, since $g^*(pdq) = d(r\circ f)$ we 
have $g^*\omega = d(g^*(pdq)) = 0$. 
Furthermore, if there is a plane in $T_x\KK^n$, for some $x \in U$,  
 included in the kernel of 
$g_* : T_x\KK^n \to T_g(x)T^*\KK^n$, then $d(p_i\circ f), d(q_i\circ f), (1 \leq i \leq n)$ 
vanish on the plane, and then also $d(r\circ f) = d(p\circ f d(q\circ f))$ vanishes 
on the plane. This means 
that the plane is included in the kernel of 
$f_* : T_x\KK^n \to T_{f(x)}\KK^{2n+1}$. Therefore if $f$ 
is of corank $\leq 1$, then $g$ is necessarily of corank $\leq 1$. 

Now, by \cite{Camb} Theorem 2, $g$ is approximated by an isotropic mapping 
$\widetilde{g} : U \to T^*\KK^n$ of corank $\leq 1$ such that, at any point $x \in 
U$, the germ of $\widetilde{g}$ at $x$ is an (symplectic) open Whitney umbrella. 
Then there exist a symplectomorphism $\kappa : (T^*\KK^n, \widetilde{g}(x)) \to 
(T^*\KK^n, 0)$ and a diffeomorphism $\sigma : (\KK^n, x) \to (\KK^n, 0)$ 
such that $\kappa\circ \widetilde{g}\circ\sigma^{-1} : (\KK^n, 0) \to (T^*\KK^n, 0)$ 
coincides with $(p\circ f_{n, k}, q\circ f_{n, k})$ in \S \ref{Main results.}. 

Let $e: (\KK^n, x) \to \KK$ be a generating function of $\widetilde{g}_x$, $\widetilde{g}_x^*(pdq) = de$. 
Remark that two generating functions $e, e'$ differ by just the addition of a constant function. 
Then $(\widetilde{g}, e) : (\KK^n, x) 
\to \KK^{2n+1}$ is an integral map-germ and it is 
contact equivalent to $f_{n, k}$ by the contactomorphism $(\sigma, \tau)$, 
$\tau(p, q, r) = (\kappa(p, q), r + c)$ for some constant $c$. 

Since $g$ is of corank $\leq 1$, if the perturbation $\widetilde{g}$ of $g$ 
is sufficiently small, then 
we can take $e$ on $U$, deleting $U$ smaller if necessary. 
Then $(\widetilde{g}, e)$ is an integral perturbation 
of $f$ on $U$, which we can take near $f$ arbitrarily. 
Since the original germ $f$ is contact stable, 
it is contact equivalent to some $(\widetilde{g}, e) : (\KK^n, x) \to \KK^{2n+1}$. 
Thus $f$ is an open Whitney umbrella in the sense of \S \ref{Main results.}. 
\enP

\noindent {\it Proof of Proposition \ref{contact stable}:} 

(cs) $\Rightarrow$ (owu) is already proved in Lemma \ref{stableOWU}.

(owu) $\Rightarrow$ (ics): 
Note that 
the infinitesimal contact stability is invariant under contactomorphisms. 
Let $f_{n,k}$ be the normal form of 
an open Whitney umbrella. Then the corresponding 
isotropic map-germ $\Pi\circ f_{n, k}$ is an open Whitney umbrella 
as an isotropic map-germs. 
Then it is proved in \cite{Inv} that $\Pi\circ f_{n, k}$ is 
symplectically stable. Then 
by Lemma \ref{symp-cont}, we see $f_{n, k}$ 
is infinitesimally contact stable. 

(owu) $\Rightarrow$ (ca) : It follows from Lemma \ref{prop OWU}. 

(ca) $\Rightarrow$ (ics) : It follows from Corollary \ref{prop codim 2}. 

(ics) $\Rightarrow$ (hcs) :  
The condition (ics) is equivalent to that 
$S_f = 0$, which is equivalent to that 
$S_{F/J} = m_JS_{F/J} = m_{W\times J}S_{F/J}$, by 
Lemma \ref{Quotient}. 
By Corollary \ref{Fin}, $S_{F/J}$ is a finite 
${\mathcal E}_{W\times J}$-module. 
So by Nakayama's lemma, we see 
$S_{F/J} = 0$. Therefore any integral deformation of $f$ 
is trivialised with respect to contactomorphisms. 
Thus we have (hcs).  

(hcs) $\Rightarrow$ (ct): 
This follows from Lemma \ref{cont-trans}. 

(ct) $\Rightarrow$ (owu) : It is 
proved in Proposition \ref{owu,ct}. 

Thus we see conditions 
(owu), (ca), (ics), (hcs) and (ct) are equivalent to each other. 

(ct) $\Rightarrow$ (cs) : 
If $j^{r}f$ is transversal to the contactomorphism class 
of $j^rf(x_0)$ for $r \geq \dfrac{n}{2} + 1$, 
then, for any slight perturbation $f'$ of $f$, there exists a point $x_0'$ near $x_0$ 
such that $j^rf'$ intersects to the contactomorphism class 
of $j^rf(x_0)$ at $x_0'$. Since $f$ is an open Whitney umbrella,  
$f$ is $r$-determined by contactomorphisms. 
Therefore we see 
$f' : (N, x_0') \to W$ is 
contactomorphic to $f : (N, x_0) \to W$. Therefore $f$ is contact stable. 

(cs) $\Rightarrow$ (ct) : 
Take a representative 
$f : U \to W$ of $f$. 
Then $f$ is approximated by an integral mapping $f' : U \to W$ 
such that $j^rf' : U \to J^r_I(N, W)$ is transverse to 
the contactomorphism-orbit $[j^rf(x_0)]$. 
Since $f$ is contact stable, there exists $x_0' \in U$ such that 
$f' : (N, {x_0'}) \to W$ and $f : (N, x_0) \to W$ are contactomorphism. 
Then $j^rf'$ is transverse to $[j^rf(x_0)] = [j^rf'(x_0')]$ at $x_0'$, and therefore 
$j^rf$ is transverse to $[j^rf(x_0)]$ at $x_0$.

\QED

\ 

Based on 
Proposition \ref{contact stable}, 
now we prove the main result Theorem \ref{Stability}.  

\ 

\noindent{\it Proof of Theorem {\ref{Stability}}} : 

First we show the equivalence of (hs) and (is). 

(hs) $\Rightarrow$ (is): 
Let $v \in VI_f$. 
Then there exists an integral deformation 
$(f_t)$ of $f$ with 
$(df_t/dt)\vert_{t=0} = v$. Since $f$ is homotopically Legendre stable, 
$f_t$ is trivialised under Legendre equivalence: 
$f_t = \tau_t^{-1}\circ f \circ\sigma_t$. 
Differentiating both sides by $t$ and setting $t = 0$, we have 
$v = (df_t/dt)\vert_{t=0} = tf(\xi) + wf(\eta)$, 
for some $\xi \in V_N, \eta \in VL_{W}$. Thus we have (is).  

(is) $\Rightarrow$ (hs): 
Since $f$ is infinitesimally Legendre stable, $f$ 
is infinitesimally contact stable. 
So $f$ is an open Whitney umbrella and thus $f$ is finite. 
Therefore $R_f$ is a finite ${\mathcal E}_N$-module. 
Then 
$VI_f/wf(VL_{W})$ is a finite ${\mathcal E}_Z$-module. 
Let $f_t$ be an integral deformation of $f$. Set $F = (f_t, t)$. 
Then 
$VI_{F/J}/(wF/J)(VL_{W})$ is also a finite ${\mathcal E}_{Z\times J}$-module. 
Thus, by Nakayama's lemma, we have 
$VI_{F/J}/((wF/J)(VL_{W}) + (tF/J)(V_N)) = 0$, similarly to the 
proof of Proposition \ref{contact stable}.  
Therefore $f$ is homotopically Legendre stable. 

Second we show (hs) ($\Leftrightarrow$ (is)) $\Rightarrow$ 
(t$_r$) $\Rightarrow$ (a$_r''$) $\Rightarrow$ (a$'$) $\Rightarrow$ 
(a) $\Rightarrow$ (is). (Therefore these conditions are equivalent to each other). 

(hs) $\Rightarrow$ (t$_r$): 
It is clear since the condition (t$_r$) is equivalent to that 
$$
VI_f = tf(V_N) + wf(VL_{W}) + VI_f^r, 
$$
by Lemma \ref{transversality cond}. 

(t$_r$) $\Rightarrow$ (a$''_r$): 
Taking generating functions of both sides of the equality 
$$
VI_f = tf(V_N) + wf(VL_{W}) + VI_f^r, 
$$
we have 
$$
R_f = (\pi\circ f)^*{\mathcal E}_Z 
+ \sum_{i=1}^n (\pi\circ f)^*{\mathcal E}_Z (p_i\circ f) 
+ R_f \cap m_N^{r+2}. 
$$
Remarking $R_f = f^*{\mathcal E}_{W}$, we have (a$''_r$). 

(a$''_r$) $\Rightarrow$ (a$'$): 
Since $R_f \cap m_N^{r+2} \subset m_W^{n+3}R_f$, (a$''_r$) implies that 
$R_f/(m_ZR_f + m_W^{n+3}R_f)$ is generated by 
$1, p_1\circ f, \dots, p_n\circ f$ over $\KK$. Then we have 
$m_W^{n+1}R_f \subset m_ZR_f + m_W^{n+3}R_f$, therefore, by Nakayama's lemma, 
$m_W^{n+1}R_f \subset m_ZR_f$. Then $m_ZR_f + m_W^{n+3}R_f = m_ZR_f$, 
so we have that $R_f/m_ZR_f$ is generated by 
$1, p_1\circ f, \dots, p_n\circ f$ over $\KK$, namely, the condition (a$'$).

(a$'$) $\Rightarrow$ (a): 
By the assumption, and by the Malgrange's preparation theorem 
of differentiable algebras (\cite{Malgrange}), 
we see $R_f$ is generated by $1$ and $p_1\circ f, \dots, p_n\circ f$ 
over ${\mathcal E}_Z$. 

(a) $\Rightarrow$ (is): Since $f$ is an open Whitney umbrella, 
we have $R_f = f^*{\mathcal E}_{W}$, and 
we have an ${\mathcal E}_Z$-isomorphism 
$VI_f/(tf(V_N) + wf(VL_{W})) 
\cong R_f/f^*({\mathcal E}_Z + \sum_{i=1}^n {\mathcal E}_Z\cdot p_i) 
$. 
Thus we see (a) implies (is). 

Lastly we show (s) $\Rightarrow$ (t$_r$)($\Leftrightarrow$ (is)) 
$\Rightarrow$ (s). 

(s) $\Rightarrow$ (t$_r$): 
Take a representative 
$f : U \to W$ of $f$. 
Then $f$ is approximated by an integral mapping $f' : U \to W$ 
such that $j^rf' : U \to J^r_I(N, W)$ is transverse to 
the Legendre orbit $[j^rf(x_0)]$. 
Since $f$ is Legendre stable, there exists $x_0' \in U$ such that 
$f' : (N, {x_0'}) \to W$ and $f : (N, x_0) \to W$ are Legendre equivalent. 
Then $j^rf'$ is transverse to $[j^rf(x_0)] = [j^rf'(x_0')]$ at $x_0'$, and therefore 
$j^rf$ is transverse to $[j^rf(x_0)]$ at $x_0$.

(t$_r$) \& (is) $\Rightarrow$ (s): 
If $j^rf$ is transverse to $[j^rf(x_0)]$ at $x_0$, then 
there exists a neighborhood $\Omega  
\subseteq C^{\infty}_I(N, W)^1$ of 
an integral representative $f : U 
\to W$
such that, for any $f' \in \Omega$, $j^rf'$ is 
transverse to $[j^rf(x_0)]$ at a point $x_0' \in U$. 
Since $j^rf'(x_0') \in [j^rf(x_0)]$, there exists an integral 
map-germ $f'' : (X, x_0) \to W$ which is Legendre equivalent to 
$f'_{x_0'}$ with respect to $\pi$ and $j^rf''(x_0) = j^rf(x_0)$. 
On the other hand, since $f$ is infinitesimally Legendre stable, 
$f$ is Legendre $r$-determined (Lemma \ref{r-determined}(4)).  
Therefore $(f'', \pi)$ and $(f, \pi)$ are Legendre equivalent. Thus 
$(f'_{x_0'}, \pi)$ and $(f, \pi)$ are Legendre equivalent, and 
$f$ is Legendre stable. 

Thus we have proved Theorem \ref{Stability}. 

\section{Contact and Legendre versalities.} 
\label{versality}

The basic singularity theory originated by H. Whitney, 
R. Thom, J. Mather, J. Martinet, C.T.C. Wall and other peoples, 
are, in particular, unified into the theory of 
geometric subgroups of ${\mathcal A}$ or ${\mathcal K}$ 
due to J. Damon \cite{Damon0}\cite{Damon1}\cite{Damon2}. 
Naturally 
we try to apply the theory of 
differentiable mappings to our situation. 
The Damon's theory guarantees the unfolding theorem (the versality theorem) 
and the determinacy theorem for a subgroup ${\GG}$ 
of ${\mathcal{A}}$ or 
${\mathcal{K}}$ acting on a linear subspace ${\FF}$ of map-germs 
${\mathcal E}(n, p) = \{ f : (\KK^n, 0) \to (\KK^p, 0), \ C^\infty\}$, 
provided that $\GG$ and $\FF$ together with their \lq\lq unfolding spaces" 
$\GG_{un}$, $\FF_{un}$ 
satisfy several required conditions. 

However our space 
$$
{\mathcal F} = \{ f : (\KK^n, 0) \to (\KK^{2n+1}, 0) \mid f \ {\mbox{\rm is integral of 
corank at most one.}} \}
$$ 
is not linear. 
Therefore, we can not apply directly the ordinary theory to our case. 

There are two possibilities to overcome this difficulty. 

One is the reduction to the linear situations case by case. For example, 
the method of generating families, due to H\" ormander and Arnold, 
is successful for the study of singularities of Lagrange and Legendre 
immersions. Moreover the linear theory successfully is applied to 
certain non-linear spaces such as spaces of solutions to non-linear 
partial differential equations, e.g. Hamilton-Jacobi equations, 
non-linear diffusions, and so on \cite{Damon4}. 
Note that in \cite{Damon4} also results on the finite determinacy 
are given. 

Another is the modifying of the original theory itself. 
It is useful to find a system of axioms which guarantees 
the versality theorem in non-linear cases, 
since then it is sufficient to just check the system of axioms. 
Then we observe,  under an additional axiom, 
that the same proof of the versality theorem in the original theory works well 
for the generalisation (Theorem 9.3 of \cite{Damon1}). 
Thus we give a direct generalisation of  Damon's theory to 
the non-linear situations. 
The generalisation is well-applied 
at least for Lagrange and Legendre 
singular immersions (isotropic and integral mappings) of corank $\leq 1$. 

We recall the theory on versal unfoldings: 
Groups of diffeomorphisms and spaces of mappings involve in the theory. 
Moreover we treat groups of unfoldings of 
diffeomorphisms and spaces of unfoldings of mappings.  

\ 

Let $\KK = \R$ or $\C$,  $C^\infty$ or holomorphic. 
Take a group of diffeomorphisms 
${\mathcal G} \subset \tilde{\mathcal K}$ where 
$$
\begin{array}{lr}
\tilde{\mathcal K} := \{ h : \KK^n\times\KK^p \to \KK^n\times\KK^p 
\mid & {\mbox{\rm fiber-preserving diffeomorphism-germs, }} \\ 
 & {\mbox{\rm w.r.t. the 
fibration \ }} \KK^n\times\KK^p 
\to \KK^n\}
\end{array}
$$
and a space of mappings 
${\mathcal F} \subset {\mathcal E} := \{ f : \KK^n \to \KK^p \mid 
{\mbox{\rm map-germs}} \}$. 
Let $f \in {\mathcal E}$ and $h \in \tilde{\mathcal K}$. Then, 
$h({\mbox{\rm graph}}(f)) = {\mbox{\rm graph}}(h(f))$ for the unique 
$h(f) \in {\mathcal E}$. We assume, for $f \in {\mathcal F}$ and $h \in {\mathcal G}, 
h(f) \in {\mathcal F}$. 

Furthermore we assume there are given a group of unfoldings of 
diffeomorphisms ${\mathcal G}_{un}(r) \subseteq \tilde{\mathcal K}_{un}(r)$ 
which acts on a space of unfoldings of mappings 
${\mathcal F}_{un}(r) \subseteq {\mathcal E}_{un}(r), r = 0, 1, 2, \dots$, with 
${\mathcal G}_{un}(0) = {\mathcal G}, {\mathcal F}_{un}(0) = {\mathcal F}$.  
Here  
$\tilde{\mathcal K}_{un}(r)$ is the space of 
$r$-parameter unfoldings of elements in $\tilde{\mathcal K}$, 
and 
${\mathcal E}_{un}(r)$ is the space of $r$-parameter unfoldings of elements 
in ${\mathcal E}$. 

First we assume ${\mathcal F}_{un}(r) \subset {\mathcal E}_{un}(r)$ is a linear subspace, 
relatively to the ordinary vector-space structure on ${\mathcal E}_{un}(r)$, $r = 0, 1, 2, 
\dots.$ 
Then, according to J. Damon, 
$({\mathcal G}, {\mathcal G}_{un}; {\mathcal F}, {\mathcal F}_{un})$ 
is called a {\it geometric subgroups and subspaces} if 
it satisfies the axioms: \\ 
(1) Naturality, (2) Tangent space structure, (3) Exponential map, 
(4) Filtration. \\ 
Then Damon has shown that 
axioms (1), (2), and (3) implies ${\mathcal G}$-versality theorem in ${\mathcal F}$, and 
axioms (1), (2), (3), and (4) implies  ${\mathcal G}$-determinacy theorem in ${\mathcal F}$. See \cite{Damon1}. See also \cite{Damon2}. 

Remark that $({\mathcal A}, {\mathcal A}_{un}; {\mathcal E}, {\mathcal E}_{un}), 
({\mathcal K}, {\mathcal K}_{un}; {\mathcal E}, {\mathcal E}_{un})$ and 
$(\widetilde{\mathcal K}, \widetilde{\mathcal K}_{un}; {\mathcal E}, {\mathcal E}_{un})$ are geometric. 
Also equivariant diffeomorphisms and mappings provide examples 
of geometric subgroups and subspaces. 
Damon and all predecessors formulated the theory explicitly for 
linear spaces (of non-linear mappings). 
However naturally the theory works also 
for non-linear mapping 
spaces (non-linear spaces of non-linear mappings) 
in ${\mathcal E}$ where geometric subgroups of $\widetilde{\mathcal K}$ 
act.

Consider non-linear ${\mathcal F} \subset {\mathcal E}$ with non-linear 
${\mathcal F}_{un} \subset {\mathcal E}_{un}$ with the action of 
a ${\mathcal G}_{un} \subset \tilde{\mathcal K}$. 

Consider the restriction ${\mathcal F}_{un}(r+s) 
\to {\mathcal F}_{un}(r)$ to the first $r$-parameters 
(resp. ${\mathcal F}_{un}(r+s) 
\to {\mathcal F}_{un}(s)$ to the last $s$-parameters)  
and the restriction ${\mathcal F}_{un}(r) \to {\mathcal F}$ 
(resp. ${\mathcal F}_{un}(s) \to {\mathcal F}$) 
to the origin of the parameter space: 
$$
\begin{array}{ccc}
{\mathcal F}_{un}(r+s) & \mapupright{\rest} &  {\mathcal F}_{un}(r) \\
 & & \\
\mapdown{\rest} &  & \mapdown{\rest} \\
 & & \\
{\mathcal F}_{un}(s)  & \mapdownright{\rest} & {\mathcal F}
\end{array}
$$
Now we pose: 

(3') Extension axiom: The natural mapping 
$$
{\mathcal F}_{un}(r+s) \to 
{\mathcal F}_{un}(r)\times_{\mathcal F}{\mathcal F}_{un}(s)
$$ 
to the fiber product 
is surjective, for any non-negative integers $r, s$. 

The axiom (3') states that a deformation of a $f \in {\mathcal F}$ over 
$\KK^r\times\{0\} \cup \{0\}\times\KK^s$ 
extends to a deformation over $\KK^{r+s}$ near $0$. 

By the same argument as in \cite{Damon1}, we can show that, 
if $({\mathcal G}, {\mathcal F})$ satisfies 
axioms (1), (2), (3) and (3') implies that ${\mathcal G}$-versality theorem in ${\mathcal F}$ 
holds, namely we have  the infinitesimal characterization, the 
existence and uniqueness 
of ${\mathcal G}$-versal unfoldings in ${\mathcal F}$. 

\ 

Here we show how to modify the conditions (1), (2) and (3) in 
\cite{Damon1} pp. 40--42. 

For the naturality, we need no change: 
(1) For any $\Sigma \in \GG_{\un}(r), F 
\in \FF_{\un}(r)$ and for any map-germ 
$\varphi : (\KK^s, 0) \to (\KK^r, 0), 
\varphi \in \EE(s, r)$, we have 
$\varphi^*\Sigma \in \GG_{\un}(s)$ and 
$\varphi^*F \in \FF_{\un}(s)$. 

For the tangent space structure, since 
we can not suppose $T\FF_{\un} = \FF_{\un}$ 
in the non-linear case, we have to modify the condition slightly: 
First we define the extended tangent spaces 
$T_1\GG_{\un, e}(r)$ and 
$T_F\FF_{\un, e}(r)$ from $\GG_{\un}(r+1), \FF_{\un}(r+1)$ 
in the same way as \cite{Damon1}, p.40. Then, 
(2) There exists an adequately ordered system 
of differentiable-analytic (DA) algebras $\{ R_\alpha\}$ 
in $\EE_{n+p}$ such that 
$T_1\GG_{\un, e}$ (resp. $T_F\FF_{\un, e}$) is 
a finitely generated $\{R_{\alpha, \lambda} \}$-module 
containing $T_1\GG_{un}$ (resp. $T_F\FF_{\un}$) as a finitely 
generated submodules, $\lambda$ indicating the 
parameter and that, for the extended orbit mapping 
$\alpha_F : \GG_{\un, e} \to \FF_{\un, e}$,  
the differential mapping 
$(\alpha_F)_* : T_1\GG_{\un, e} \to T_F\FF_{\un, e}$ 
is an $\{R_{\alpha, \lambda} \}$-module homomorphism. 
The finiteness condition is required only when 
$f = F\vert_{\KK^n\times 0}$ satisfies $\dim_{\KK}(T_f{\mathcal F}/T\GG_e\dot f) 
< \infty$. 
Moreover there exist isomorphisms 
$T_1\GG_{\un, e}/m_\lambda T_1\GG_{\un, e} \cong T_1\GG_e$, 
$T_F\FF_{\un, e}/m_\lambda T_F\FF_{\un, e} \cong T_f\FF_e$, 
as 
$\{R_\alpha\}$-modules, 
and that 
$\{m_\alpha\}T_1\GG_e \subset T_1\GG$, 
and 
$\{m_\alpha\}T_f\FF_e \subset T_f\FF$. 
About the generalities on DA-algebras see \cite{Damon1}\cite{Damon2}. 

For the exponential property (3) we need no change.

\bee
{\rm 
Let ${\mathcal I}$ be a differential system (= an ideal of differential forms that is 
$d$-closed) on $\KK^p$. 
Set ${\mathcal F} := \{ f : (\KK^n, 0) \to (\KK^p, 0) \mid f^*{\mathcal I} = 0\}$, 
the set of integral map-germs, and 
${\mathcal G} := \{(\sigma, \tau) \in {\mathcal A} \mid \tau^*{\mathcal I} = {\mathcal I}\}$, 
the group of ${\mathcal A}$-equivalences preserving ${\mathcal I}$. 
Then ${\mathcal F} \subset {\mathcal E}$ is ${\mathcal G}$-invariant. 
Moreover we set 
$
{\mathcal F}_{un}(r)
$ as the space of $r$-parameter unfoldings  
$$
F = (f_\lambda, \lambda) : 
(\KK^n\times\KK^r, (0,0)) \to (\KK^p\times \KK^r, (0, 0)) 
$$
satisfying 
$f_{\lambda}^*{\mathcal I} = 0  (\lambda \in (\KK^r, 0))$, 
and set 
$$
{\mathcal G}_{un}(r) := \{(\sigma_{\lambda}, \tau_{\lambda}, \lambda) 
\in {\mathcal A}_{un}(r) 
 \mid \tau_{\lambda}^*{\mathcal I} = {\mathcal I},  (\lambda \in (\KK^r, 0)) \}. 
$$
Then ${\mathcal F}_{un}(r) \subset {\mathcal E}_{un}(r)$ is ${\mathcal G}_{un}(r)$-invariant. 
Remark that ${\mathcal F}$ and ${\mathcal F}_{un}$ 
are in general non-linear. 

In particular we apply the above general theory to the singularity 
theory of integral mappings. 

Set 
$$
{\mathcal F} = \{ f : (\KK^n, 0) \to (\KK^{2n+1}, 0) \mid f \ {\mbox{\rm is integral of 
corank at most one.}} \}
$$ 
and ${\mathcal G} = \{ (\sigma, \tau)\}$ 
the group of contactomorphisms acting on ${\mathcal F}$. 
We set 
${\mathcal F}_{un}(r)$ as the space of $r$-parameter integral unfoldings 
of integral germs in ${\mathcal F}$, and 
${\mathcal G}_{un}(r)$ as the group of $r$-parameter unfoldings 
of contactomorphisms in ${\mathcal G}$. 
Then we have $T_f{\mathcal F}_e = VI_f, T_f{\mathcal F} = VI_f^0$ 
and 
$T_1{\mathcal G}_e = V_n\oplus VH_{2n}, 
T_1{\mathcal G} = m_nV_n\oplus(m_{2n}^{(2)}*VH_{2n})$. 

Then $T_f{\mathcal F}_e$ is an ${\mathcal E}_{2n+1}$-module. 
If $f$ is finite, then $T_f{\mathcal F}_e$ is a finite ${\mathcal E}_{2n+1}$-module. 

We consider the system of algebras: ${\mathcal E}_{2n+1} \to 
{\mathcal E}_{2n+1}$ with the identity connection homomorphism. 
Note that $T_1{\mathcal G}_e = V_n\oplus VH_{2n}$ is an 
$\{ {\mathcal E}_{2n+1}, {\mathcal E}_{2n+1}\}$-module 
by 
$$
H*\xi := (f^*H)\cdot\xi, \ H*X_K := X_{HK}, 
$$
where $X_H$ is the contact Hamilton vector field with Hamilton function $H$
(\S \ref{Contact Hamilton vector fields.}). 
Moreover $T_f{\mathcal F}_e$ is an 
$\{ {\mathcal E}_{2n+1}, {\mathcal E}_{2n+1}\}$-module by 
$$
\{ H, K \}* v := H*v, 
$$ 
using the multiplication given in Proposition \ref{module}. 
(Here $\{ H, K\}$ does not mean the Poisson product, but just comes 
from the notion on algebra-systems used in \cite{Damon1} \S 6). 

Also for unfolding spaces, module structures are defined as in Corollary 
\ref{Fin}. 

\ 

For the Legendre versality, we set ${\mathcal G}$ 
the group of Legendre equivalences $\{ (\sigma, \tau)\}$ 
for the Legendre fibration $\pi : W = (\KK^{2n+1}, 0) \to Z = (\KK^{n+1}, 0)$. 
We read as 
$$
VI_f = T_{f}{\mathcal F}_e {\mbox{\rm \ and \ }} 
tf(V_N) + wf(VL_W) = T{\mathcal G}_e\cdot f. 
$$
The system of algebra we consider is 
$${\mathcal E}_{n+1} \to  {\mathcal E}_{2n+1} \to {\mathcal E}_{2n+1}$$ 
with the connection homomorphisms $\pi^*$ and the identity respectively. 

In both case we can check the axioms (1), (2) and (3). 

}
\ene

\ 

Now, based on the above general scheme due to Damon after the modification, we give 
alternative proof of Theorem \ref{Versality}

\ 

\noindent{\it Proof of Theorems \ref{contact Versality} and \ref{Versality}:} 
The axioms (1), (2), (3) are easily checked. 
We need to show the extension axiom (3') is satisfied.  
for the category of integral unfoldings of an integral map-germs 
of corank at most one,  in order to apply the Mather-Damon's machine 
to our situation. Note that the geometric group does not involve into 
the axiom (3'). 

Let $f : (N, x_0) \to (W, w_0)$ be an integral map-germ of corank $\leq 1$. 
Let $F : (N \times \KK^r, (x_0, 0)) 
\to (W, w_0)$ and 
$F' : (N \times \KK^s, (x_0, 0)) 
\to (W, w_0)$ be integral deformations of $f$. 
We may set $(q, p_n)\circ f = (x_1, \dots, x_{n-1}, u, v)$, for 
a function-germs $u = u(x', t), v = v(x', t), x' = (x_1, \dots, x_{n-1}), 
t = x_n$, 
after a contactomorphism. 
Set 
$$
(q, p_n)\circ F = ((q\circ F)(x', t, \lambda), (p_n\circ F)(x', t, \lambda)), 
$$
$$
(q, p_n)\circ F' = ((q\circ F')(x', t, \mu), (p_n\circ F')(x', t, \mu)). 
$$
Then there exist a coordinate change on the $q$-space 
depending on $\lambda, \mu$ such that we have 
$$
(q, p_n)\circ F = (x', U(x', t, \lambda), V(x', t, \lambda)), 
$$
$$
(q, p_n)\circ F' = (x', U'(x', t, \mu), V'(x', t, \mu)), 
$$
with $U(x', t, 0) = u, V(x', t, 0) = v, U'(x', t, 0) = u, V'(x', t, 0) = v$. 
Then we can extend 
$(q, p_n)\circ F$ and $(q, p_n)\circ F'$ to 
$H : (N \times \KK^r\times \KK^s, (x_0, 0, 0)) 
\to (W, w_0)$ 
of form 
$$
H(x', t, \lambda, \mu) = (x', \widetilde{U}(x', t, \lambda, \mu), 
\widetilde{V}(x', t, \lambda, \mu))
$$
by setting 
$$
\widetilde{U}(x', t, \lambda, \mu) := 
U(x', t, \lambda) + U'(x', t, \mu) - u(x', t), 
$$
$$
\widetilde{V}(x', t, \lambda, \mu) := 
V(x', t, \lambda) + V'(x', t, \mu) - v(x', t). 
$$
Then we define 
$F'' : (N\times \KK^r\times \KK^s, (x_0, 0, 0)) 
\to (W, w_0)$ by $(q, p_n)\circ F'' = H$, 
$(r\circ F'')(x_0) = r\circ f(x_0)$, and by 
$$
d(r\circ F'') = \sum_{i=1}^{n-1}(p_i\circ F'')dx_i + \widetilde{U}d\widetilde{V}. 
$$
Here $d$ means the exterior differential by $x_1, \dots, x_{n-1}, x_n = t$. 
The last condition means that  
$$
\dfrac{\pa r\circ F''}{\pa x_i} = p_i\circ F'' + \widetilde{U}\dfrac{\pa \widetilde{V}}{\pa x_i}, \quad \dfrac{\pa r\circ F''}{\pa t} = \widetilde{U}\dfrac{\pa \widetilde{V}}{\pa t}. 
$$
We determine $r\circ F''$ by the latter condition and 
$(r\circ F'')(x_0) = r\circ f(x_0)$. 
Then $p_i\circ F'', (1 \leq i \leq n-1)$ are determined by the former condition. 
Thus we have the extension $F''$ of both $F$ and $F'$. 

Therefore the extension axiom (3') is satisfied. Then by the general 
framework we have the proof of contact and Legendre versality theorems. 
\QED


\ 

Go-o ISHIKAWA

Department of Mathematics, Hokkaido University, 

Sapporo 060, JAPAN 

e-mail: ishikawa@math.sci.hokudai.ac.jp 

\end{document}